\documentclass[a4paper,12pt,oneside]{article}
\usepackage{main}

\title{%
Surrogate to Poincar\'e inequalities on manifolds for dimension reduction in nonlinear feature spaces}
\author{%
  A. Nouy${}^{1}$,
  A. Pasco${}^{1}$}
\date{\medskip%
  \small %
  ${}^1$ \'Ecole Centrale Nantes, Nantes Universit\'e, \\ Laboratoire de Mathématiques Jean Leray UMR CNRS 6629\\
  \texttt{alexandre.pasco1702@gmail.com;anthony.nouy@ec-nantes.fr}
}

\begin{document}
\maketitle

\begin{abstract}

We aim to approximate a continuously differentiable function $u:\mathbb{R}^d \rightarrow \mathbb{R}$ by a composition of functions $f\circ g$ where $g:\mathbb{R}^d \rightarrow \mathbb{R}^m$, $m\leq d$, and $f : \mathbb{R}^m \rightarrow \mathbb{R}$ are built in a two stage procedure.
For a fixed $g$, we build $f$ using classical regression methods, involving evaluations of $u$.
Recent works proposed to build a nonlinear $g$ by minimizing a loss function $\mathcal{J}(g)$ derived from Poincar\'e inequalities on manifolds, involving evaluations of the gradient of $u$.
A problem is that minimizing $\mathcal{J}$ may be a challenging task.
Hence in this work, we introduce new convex surrogates to $\mathcal{J}$.
Leveraging concentration inequalities, we provide suboptimality results for a class of functions $g$, including polynomials, and a wide class of input probability measures.
We investigate performances on different benchmarks for various training sample sizes.
We show that our approach outperforms standard iterative methods for minimizing the training Poincar\'e inequality based loss, often resulting in better approximation errors, especially for small training sets and $m=1$.
\end{abstract}



\section{Introduction}
\label{sec:Introduction}

We consider the problem of approximating a  function $u : \bbR^d \to \bbR$ 
by a composition $f \circ g$ with $g : \bbR^d \to \bbR^m$  a feature map which associates to $\bfx\in \bbR^d$  a set of $m$ real features $g(\bfx) = (g_1(\bfx), \hdots, g_m(\bfx)) \in \bbR^m$. The map $g$ is selected from some  model class $\calG_m$. 
We assume that $u$ is in $L^2(\bbR^d , \mu)$, where $\mu$ is some probability measure, which we equip with the natural norm $\Vert u \Vert^2 := \bbE(\vert u(\bfX) \vert^2 )$, where $\bfX $ is a random vector with distribution $\mu$. 
For a fixed function $g$, an optimal $L^2$-approximation of the form $f \circ g$ is a solution of the optimization problem   
\begin{equation}
\label{equ:def reconstruction error} 
  \min_{f}
  \Expe{|u(\bfX) - f\circ g(\bfX)|^2}:=  \calE (g),
\end{equation}
where the minimum is taken over all measurable maps $f$ from $\bbR^m$ to $\bbR$. It admits as a solution  the regression function  $f^\star (\bfz) = \bbE( u(\bfX) \vert g(\bfX)= \bfz)$. 
An optimal feature map $g$ from the class  $\calG_m$  can then be defined  by minimizing  the  error  $\calE (g)$ over $\calG_m$. 
This results in an optimal dimension reduction.

The simplest model class $\calG_m$ one could think of are linear feature maps $g(\bfx) = G^T \bfx$ for some matrix $G \in \bbR^{d \times m}$, which can be referred to as \textit{linear dimension reduction}, and the associated approximation of $u$ is then a ridge function \cite{pinkusRidgeFunctions2015}.
The statistical regression literature denoted this framework as (linear) \textit{sufficient dimension reduction} \cite{adragniSufficientDimensionReduction2009}, in which the main idea is to find a small $m$-dimensional subspace spanned by the columns of $G \in\bbR^{d\times m}$, such that the random variables $u(\bfX)$ and $\bfX | G^T \bfX$ are independent.
This subspace is generally built using methods based on conditional variance, such as sliced inverse regression (SIR) \cite{liSlicedInverseRegression1991}, sliced average variance estimation (SAVE) \cite{cookSlicedInverseRegression1991} or directional regression (DR) \cite{liDirectionalRegressionDimension2007}.
A main issue of such methods is that they do not provide certification on the approximation error $\calE(g)$.
On the other hand, the so-called active subspace method \cite{constantineActiveSubspaceMethods2014,zahmGradientBasedDimensionReduction2020}, which is a gradient-based method closely related to the present work, provides such a certification by leveraging probabilistic Poincar\'e inequalities.
The latter yields an upper bound of the form
\begin{equation}
\label{equ:poincare based bound intro}
  \calE(g) \leq C(\bfX|\calG_m) \calJ(g)
\end{equation}
where $C(\bfX|\calG_m)$ is a constant depending only on the law of $\bfX$ and the function class $\calG_m$ and where $\calJ(g) = \Expe{\|\Piperp_{\spanv{G}} \nabla u(\bfX)\|_2^2}$.
Here, for any matrix $M \in \bbR^{d_1\times d_2}$ and integers $d_1,d_2>0$, $\spanv{M}$ is the column span of $M$ and $\Piperp_{\spanv{M}} := \Pi_{\spanv{M}^{\perp}} \in \bbR^{d_1\times d_1}$ is the orthogonal projector onto $\spanv{M}^{\perp}$, which is the orthogonal complement to $\spanv{M}$.
The constant $C(\bfX|\calG_m)$ is defined via probabilistic Poincar\'e constants associated to the conditional random variables $\bfX | G^T \bfX=\bfz$ for $\bfz\in\bbR^m$.
The fact that the support of the latter lies in the $(d-m)$-dimensional affine space $\bfz + \spanv{G}^{\perp}$ allows us to derive bounds on $C(\bfX|\calG_m)$ for several classical probability measures \cite{teixeiraparenteGeneralizedBoundsActive2020}.
The other main advantage of this method is that, while minimizing $\calE$ on $\bbR^{d\times m}$ is a challenging problem, minimizing $\calJ$ on $\bbR^{d\times m}$ is rather straightforward.
Indeed, it is a quadratic functional, whose minimum is the sum of the $d-m$ smallest eigenvalues of the matrix $\Expe{\nabla u(\bfX) \nabla u(\bfX)^T} \in \bbR^{d\times d}$, and the  minimizer is a matrix whose columns are  the eigenvectors associated with the $m$ largest eigenvalues of that matrix.

Although linear dimension reduction has several theoretical and practical advantages, the main limitation is that it may not efficiently detect important nonlinear features, as it is easy to see by considering $u(\bfx) = h(\|\bfx\|^2_2)$ for some $h\in\calC^{1}$.
In this regard, there are insights for developing \textit{nonlinear dimension reduction} methods, where the model class $\calG_m$ contains nonlinear functions.
Until recently, there were mainly \textit{nonlinear sufficient dimension reduction} methods developed,
see for example \cite{yi-renyehNonlinearDimensionReduction2009,leeGeneralTheoryNonlinear2013,liNonlinearSufficientDimension2017,liDimensionReductionFunctional2022}, extending the linear sufficient dimension reduction methods (SIR, SAVE, DR).
More recently, various works introduced gradient-based nonlinear dimension reduction methods, either in a  heuristic manner \cite{zhangLearningNonlinearLevel2019,bridgesActiveManifoldsNonlinear2019,romorKernelbasedActiveSubspaces2022}, or using Poincar\'e inequalities  \cite{bigoniNonlinearDimensionReduction2022,verdiereDiffeomorphismbasedFeatureLearning2025}.  
In particular in \cite{bigoniNonlinearDimensionReduction2022} it was shown that \eqref{equ:poincare based bound intro} still holds, under important assumptions, with
\begin{equation}
  \calJ(g) = \Expe{\|\Piperp_{\spanv{\nabla g(\bfX)}} \nabla u(\bfX)\|^2_2}
  = \Expe{\|\nabla u(\bfX)\|^2_2} -  \Expe{\|\Pi_{\spanv{\nabla g(\bfX)}} \nabla u(\bfX)\|^2_2},
\end{equation}
where $\nabla g(\bfX) := (\nabla g_1(\bfX), \cdots \nabla g_m(\bfX)) = (\frac{\partial g_j}{ \partial x_i} (\bfX))_{1\leq i\leq d, 1\leq j\leq m} \in \bbR^{d\times m}$.
The constant $C(\bfX|\calG_m)$ is here defined via probabilistic Poincar\'e inequalities on manifolds, associated to the conditional random variables $\bfX | g(\bfX)=\bfz$ for $\bfz\in\bbR^m$.
Indeed, the support of the latter lies in $g^{-1}(\{\bfz\})$, which is no longer an affine space but rather a $(d-m)$-dimensional smooth submanifold of $\bbR^d$, under some assumptions on $g$. 
Note that $\Piperp_{\spanv{\nabla g(\bfX)}} \nabla u(\bfX)$ is actually the Riemannian gradient of $u_{|g^{-1}(\{g(\bfX)\})}$ at point $\bfX$, associated to the Euclidean metric in $\bbR^d$.

This nonlinear setting comes with two main issues.
The first is that the theoretical foundation of \eqref{equ:poincare based bound intro} is much   weaker than for linear features because controlling, or simply ensuring finiteness, for $C(\bfX|\calG_m)$ is a truly challenging task.
A first step in overcoming this problem has been proposed in \cite{verdiereDiffeomorphismbasedFeatureLearning2025} by considering diffeomorphisms-based feature  maps $g: \bfx\mapsto (\varphi_1(\bfx), \cdots, \varphi_m(\bfx))$, with $\varphi$ a $\calC^1$ diffeomorphism on $\bbR^d$.
The second issue is that optimizing $\calJ$ over a class $\calG_m$ of nonlinear functions can be a challenging task, although first and second order optimization algorithms are available.
This is mainly due to the nonlinearity of the mapping $g \mapsto \Piperp_{\spanv{\nabla g(\bfx)}}$.
This second issue is the motivation for the present work, in which we derive new upper bounds that yield to straightforward optimization procedures.

Finally, we emphasize the fact that the approaches described in this introduction, as well as the one developed in the present paper, are two steps procedures.
Firstly, we build the feature map $g$.
Secondly, since the best regression function $f^*$ is not available in practice, we build $f:\bbR^m \rightarrow \bbR$ by minimizing $\Expe{(u(\bfX) - f\circ g(\bfX))^2}$ over some regression class $\calF_m$.
As a result, the construction of $g$ is oblivious to the class $\calF_m$, which might be problematic in practice.
Another approach would be to optimize on $f$ and $g$ in a coupled way, as for example proposed in \cite{hokansonDataDrivenPolynomialRidge2018,lataniotisEXTENDINGCLASSICALSURROGATE2020}.
We also emphasize the fact that the approaches based on Poincar\'e inequalities require samples of $u(\bfX)$ as well as $\nabla u(\bfX)$.

\paragraph{Contributions.}

The main contribution of the present work is to introduce new surrogates in order to circumvent the challenge of optimizing $\calJ$.
Those surrogates are  much easier to optimize, and their suboptimality with respect to $\calJ$ can be somehow controlled, under some assumptions on $\calG_m$.
The first assumption that we will require in the present work is the following.
\begin{assumption}
\label{assump:full rank jacobian as}
  For all $g\in\calG_m$ it holds rank$(\nabla g(\bfX))=m$ almost surely.
\end{assumption}
This first assumption will allow us to show recovery properties of one feature when using our surrogates.
More precise suboptimality results will be obtained by further assuming that distributions associated with $\calG_m$ satisfy some \textit{deviation inequalities}.
We mainly provide results for the one feature case $m=1$, for which we introduce the following new surrogate,
\[
  \calL_1(g)
  = \Expe{\|\nabla u(\bfX)\|_2^2 \|\Piperp_{\spanv{\nabla u(\bfX)}} \nabla g(\bfX) \|_2^2},
\]
which is quadratic in $g$.
If $\calG_1$ satisfies some assumptions allowing to draw deviation inequalities on $\|\nabla g(\bfX)\|_2^2$, typically polynomials, then we show that
\[
  \calL_1(g) \lesssim \calJ(g) \lesssim \calL_1(g)^{\beta},
\]
with $\beta \in (0,1]$ depending on the rate of decay involved in the deviation inequalities.
We also show that that if $g(\bfx) = G^T \Phi(\bfx)$ for some $\Phi\in\calC^1(\bbR^d, \bbR^K)$ and some $G\in\bbR^K$ then $\calL_1(g) = G^T H G$ with
\[
  H = \Expe{
    \nabla \Phi(\bfX)^T 
    \big(
      \|\nabla u(\bfX)\|_2^2 I_d - \nabla u(\bfX) \nabla u(\bfX)^T
    \big)  
    \nabla \Phi(\bfX)
    } 
    \in \bbR^{K \times K},
\]
which means that minimizing $\calL_1$ is equivalent to finding an eigenvector associated to the smallest eigenvalue of $H$.

We also propose an extension of our approach to the case of multiple features $(m>1)$.
We first introduce for any $1\leq j\leq m$ a new surrogate quantity $\calL_{m,j}:\calG_m \rightarrow \bbR$ which is quadratic with respect to $g_j$.
We show that for polynomial features this surrogate satisfies similar suboptimality properties as in the case $m=1$, but with a deteriorated exponent $\beta$. 
We then propose a greedy algorithm to sequentially learn $g_1, \cdots, g_m$, where at iteration $j$ of the algorithm we choose $g_j$ as a minimizer of $h \mapsto \calL_{m,j}((g_1, \cdots, g_{j-1}, h, 0, \cdots, 0))$.

In practice, the feature maps $g$ obtained by minimizing our surrogates for $m=1$ or $m>1$ can be used as good starting guesses of descent algorithms for minimizing $\calJ$, or directly as feature maps for  building the regression function $f:\bbR \rightarrow \bbR$.

The rest of this paper is organized as follows.
First in \Cref{sec:Dimension reduction using Poincar\'e inequalities} we summarize the   nonlinear dimension reduction method based on Poincar\'e inequalities on manifolds introduced in \cite{bigoniNonlinearDimensionReduction2022}, and we discuss on particular choices for the model class $\calG_m$.
In \Cref{sec:Large and small deviations inequalities} we shortly summarize and adapt to our framework the deviation inequalities from \cite{fradeliziConcentrationInequalities$s$concave2009}, which we leverage to ensure theoretical guarantees on our surrogates.
Then in \Cref{sec:one feature} we introduce our new surrogate $\calL_1$ for $m=1$ and provide a general analysis.
In \Cref{sec:multiple features} we extend the method to the case $m>1$ using a greedy algorithm.
Finally, in \Cref{sec:numerical experiments} we illustrate our methods on numerical examples.

\section{Dimension reduction using Poincar\'e inequalities}
\label{sec:Dimension reduction using Poincar\'e inequalities}

In this section, we detail how one can leverage Poincar\'e inequalities to construct a feature map $g$  for the approximation of $u$ by a composition $f\circ g$, with $g$ from some model class $\calG_m$ and  $f$ some measurable map from $\bbR^m $ to $ \bbR$. 
First in \Cref{subsec:Poincare Inequality based upper bound} we present the nonlinear dimension reduction approach proposed in \cite{bigoniNonlinearDimensionReduction2022}.
Then in \Cref{subsec:Choice of the feature maps} we discuss on different choices of model classes $\calG_m$.

\subsection{Error bounds based on Poincar\'e inequalities}
\label{subsec:Poincare Inequality based upper bound}

A key observation is that for any $g\in\calC^1(\calX, \bbR^m)$ and $f \in \calC^1(\bbR^m, \bbR)$, applying the chain rule to $u = f\circ g$ yields $\nabla u(\bfx) = \nabla g(\bfx) \nabla f(g(\bfx))$, where $\nabla g(\bfx) \in \bbR^{d\times m}$ is the jacobian matrix of $g$ as defined in \Cref{sec:Introduction}.
Therefore, 
\[
\nabla u(\bfx)   \in \spanv{\nabla g(\bfx)} := \spanv{\nabla g_i(\bfx)}_{1\leq i\leq m} \subset \bbR^d.
\]
Hence, there is an incentive to build $g$ such that $\nabla u(\bfx)$ and $\spanv{\nabla{g(\bfx)}}$ are as aligned as possible.
The authors in \cite{bigoniNonlinearDimensionReduction2022} propose to minimize the following cost function,
\begin{equation}
\label{equ:def of J}
  \calJ(g) : = \Expe{\|\nabla u (\bfX)\|_2^2} - \Expe{\|\Pi_{\nabla g(\bfX)} \nabla u (\bfX)\|_2^2},
\end{equation}
where we denote $\Pi_{M} := \Pi_{\spanv{M}}$ where $\Pi_{\spanv{M}} \in \bbR^{d \times d}$ is the orthogonal projector onto the column span of $M$ for any matrix $M\in\bbR^{d\times d'}$ and any integer $d'>0$.
We similarly denote $\Piperp_{M} := \Piperp_{\spanv{M}} = \Pi_{\spanv{M}^{\perp}}$.
The cost function $  \calJ(g)$ naturally appears when using Poincar\'e inequalities and provides certified error bounds for $g \in \calG_m$ \cite{zahmGradientBasedDimensionReduction2020}.
When $\calG_m$ contains nonlinear functions, these error bounds involve Poincar\'e inequalities on Riemannian manifolds, which we define in \Cref{def:poincare inequality}  as in \cite{bigoniNonlinearDimensionReduction2022,verdiereDiffeomorphismbasedFeatureLearning2025}. 
Deriving these inequalities involve concepts from  differential geometry, and although a deep understanding of the latter is not required to read the current section, we recommend \cite[Chapter 3]{absilOptimizationAlgorithmsMatrix2008} for a comprehensive introduction to differential geometry and its applications.

\begin{definition}[Poincar\'e inequality on Riemannian submanifold of $\bbR^d$]
\label{def:poincare inequality}
  Let $\bfY$ a continuous random vector taking values in a Riemannian submanifold $\calM \subset \bbR^d$ equipped with the Euclidean metric.
  Then the Poincar\'e constant $C(\bfY) \geq 0$ is the smallest constant such that
  \begin{equation}
  \label{equ:def poincare inequality}
    \Expe{(h(\bfY) - \Expe{h(\bfY)})^2}
    \leq C(\bfY) \Expe{\|\nabla h (\bfY)\|_2^2}
  \end{equation}
  holds for any $h\in\calC^1(\calM, \bbR)$. 
  Here $\nabla h $ denotes the Riemannian gradient.
  $\bfY$ is said to satisfy Poincar\'e inequality if $C(\bfY) < +\infty$.
\end{definition}

One of the differences between \eqref{equ:def poincare inequality} and the classical Poincar\'e inequality in Euclidean spaces is that it involves the Riemannian gradient $\nabla h(\bfy)\in T_{\bfy} \calM$, where $T_{\bfy} \calM$ is the tangent space of $\calM$ at $\bfy\in\calM$.
Fortunately, the fact that $\calM$ is embedded in $\bbR^d$ simplifies some concepts involved.
Indeed,   we can identify  $T_{\bfy} \calM$  with a linear subspace of $\bbR^d$ and express the Riemannian gradient of $h\in\calC^1(\calM, \bbR)$ as 
\[
  \nabla h(\bfy)
  = \Pi_{T_{\bfy} \calM} \nabla \bar{h} (\bfy),
\]
where $\Pi_{T_{\bfy} \calM} $ is the orthogonal projection onto $T_{\bfy} \calM$, and $\bar{h}\in\calC^1(\bbR^d, \bbR)$ is a function whose restriction $h_{|\calM} $  to $\calM$ is equal to $h$.

In our setting, the manifolds involved are the level-sets of $g$, which are smooth manifolds of dimensions $d-m$, assuming that $\text{rank}(\nabla g(\bfx)) = m$ for all $\bfx\in\calX$.
The manifold-valued random vector involved is the conditional random variable $\bfY = (\bfX | g(\bfX) = \bfz)$ for some $\bfz\in g(\calX)$, which takes values in the manifold $\calM_{\bfz} = g^{-1}(\{\bfz\})$ (a level-set of $g$), and the function defined on $\calM_{\bfz}$ is $u_{|\calM_{\bfz}}$.
In this case, for any $\bfx\in\calM_{\bfz}$, it holds 
\[
T_{\bfx} \calM_{\bfz} = \text{ker} \{\nabla g(\bfx)^T\} = \spanv{\nabla g(\bfx)}^\perp
\] 
and thus 
\[
  \nabla (u_{|\calM_{\bfz}})(\bfx)
  = \nabla u(\bfx) - \Pi_{\nabla g(\bfx)} \nabla u(\bfx).
\]
Now applying \eqref{equ:def poincare inequality} to $u_{|\calM_{\bfz}}$ and taking the supremum over  $\bfz\in g(\calX)$ and  $g\in\calG_m$  of the Poincar\'e constant $C(\bfX | g(\bfX)=\bfz)$ yields \Cref{prop:poincare based bound}, which is a restatement of  \cite[Proposition 2.9]{bigoniNonlinearDimensionReduction2022}.  

\begin{proposition}
\label{prop:poincare based bound}
  Let $\calG_m \subset \calC^1(\calX, \bbR^m)$ satisfying rank$(\nabla g(\bfx)) = m$ for all $g\in\calG_m$ and all $\bfx\in\calX$, and assume that  
  \begin{equation}
  \label{equ:def of uniform C}
    C(\bfX | \calG_m):= 
    \sup_{g \in \calG_m} \sup_{z \in g(\calX)} 
    C(\bfX | g(\bfX)= \bfz) < + \infty,
  \end{equation}
  with $C(\bfX | g(\bfX)= \bfz)$ defined in \eqref{equ:def poincare inequality}.
  Then for any $u\in\calC^1(\calX, \bbR)$ and $g\in\calG_m$ it holds 
  \begin{equation}
  \label{equ:poincare based bound}
    \calE(g) \leq C(\bfX | \calG_m) \calJ(g),
  \end{equation}
  with $\calE$ defined in \eqref{equ:def reconstruction error} and  $\calJ$ defined in \eqref{equ:def of J}.
\end{proposition}

The bound from \Cref{prop:poincare based bound} is the theoretical justification behind the idea of minimizing $\calJ$ in order to build the nonlinear feature map $g$.
However, this strategy might be far from optimal since we are minimizing an upper bound of the true reconstruction error $\calE$, especially if the constant $C(\bfX | \calG_m)$ is large.
Still, the numerical results in \cite{bigoniNonlinearDimensionReduction2022,verdiereDiffeomorphismbasedFeatureLearning2025} show good performances of such Poincar\'e inequalities based approach on various practical examples.

It is important to note that $\calJ$ satisfies some invariance properties.
Indeed,  for any diffeomorphism $\psi : \bbR^m \rightarrow \bbR^m$, by the chain rule and invertibility of $\nabla \psi(\bfx)$, it holds  that $\spanv{\nabla (\psi \circ g(\bfx))} = \spanv{\nabla g(\bfx)}$ and thus 
\begin{equation}
\label{equ:invariance J diffeo}
  \calJ(\psi \circ g) = \calJ(g).
\end{equation}
This defines an equivalence relation on $\calG_m$, which allows in practice to restrict the set $\calG_m$. In particular, since $\calJ(g/\|g\|) = \calJ(g)$ for any $g\neq 0$, 
one can impose that $\calG_m$ is contained in the unit sphere of some normed space. 
Also, up to a suitable  affine transformation $\psi$, one can consider that $\Expe{g(\bfX)} = 0$ and $\Expe{g(\bfX) g(\bfX)^T} = I_m$, which is imposed by the authors in \cite{bigoniNonlinearDimensionReduction2022} for improving numerical stability of algorithms. 
Note however that although $\calJ(\psi \circ g) = \calJ(g)$, approximating $\bfz \mapsto \Expe{u(\bfX)| g(\bfX) = \bfz}$ in $\calF$ can be more challenging than approximating $\bfz \mapsto \Expe{u(\bfX)| \psi \circ g(\bfX) = \bfz}$ in $\calF$.

Finally, it is worth mentioning that the method detailed in the current section can be extended naturally to vector valued functions, as detailed in \cite{verdiereDiffeomorphismbasedFeatureLearning2025}.

\subsection{Choice of  feature maps}
\label{subsec:Choice of the feature maps}

In this section we discuss classes of feature maps that have been   proposed in the framework described in \Cref{subsec:Poincare Inequality based upper bound}.

\subsubsection{Linear feature maps}

Linear features correspond to the linear space  
\begin{equation}
\label{equ:Gm linear features}
  \calG_m = \{ g:\bfx \mapsto G^T \bfx : G\in\bbR^{d\times m}\}.
\end{equation}
It corresponds to the so-called \textit{active subspace method}  \cite{constantineActiveSubspaceMethods2014,zahmGradientBasedDimensionReduction2020}. It has  been used successfully in various practical applications. 

This setting has some advantages. 
First, there exist known bounds on $ C(\bfX | \calG_m)$ for some standard probability distributions for $\bfX$.
For example, if $\bfX \sim \calN(0,I_d)$ then one can show that $C(\bfX | \calG_m)=1$.
Secondly, a minimizer   of $\calJ$  can be obtained by computing the leading eigenvectors of the matrix $\Expe{\nabla u(\bfX) \nabla u(\bfX)^T} \in \bbR^{d\times d}$.
However, choosing linear features restricts the class of functions $u$ that can be well approximated. As an example, the function $\bfx \mapsto \sin(\|\bfx\|_2^2)$ cannot be well approximated using a few linear features. 
 
\subsubsection{Feature maps from a linear space}

Feature maps can be chosen from some finite dimensional linear space
\begin{equation}
\label{equ:Gm vector space}
  \calG_m = \{ g:\bfx \mapsto G^T \Phi(\bfx) :  G \in \bbR^{K\times m}\}.
\end{equation}
Here $\Phi \in \calC^1(\bbR^d, \bbR^K)$, with $K\geq d$, is some fixed (or adaptively selected) feature map. 
The coordinates $g_i$ of $g\in   \calG_m$ are in the 
$K$-dimensional linear space $\spanv{\Phi_1(\bfx), \hdots, \Phi_K(\bfx)}$.  
This has been considered in 
\cite{bigoniNonlinearDimensionReduction2022,romorKernelbasedActiveSubspaces2022}.
It has the advantage of being potentially able to approximate functions with fewer features than with \eqref{equ:Gm linear features}, while keeping a  simple parameterization of $\calG_m$.
Moreover, we can compute analytically first and second order derivatives of $G \mapsto \calJ(G^T \Phi)$ and exploit the invariance property \eqref{equ:invariance J diffeo} to add constrains to the minimization procedure.
In particular $G \mapsto \calJ(G^T \Phi)$ actually only depends on $\spanv{G} \subset \bbR^K$.
Thus, one can for example consider instead the compact set 
\begin{equation}
\label{equ:Gm vector space sphere}
  \calG_m = 
  \Big\{ 
    g: \bfx \mapsto G^T \Phi(\bfx) :  G \in \bbR^{K \times m},
     G^T R G = I_m
  \Big\},
\end{equation}
with $R \in \bbR^{K\times K}$ a symmetric positive definite matrix.
We discuss on how to choose $R$ in \Cref{subsec:minimizing the surrogate}.
The set being identified with a Stiefel matrix manifold, see for example \cite[Section 3.3.2]{absilOptimizationAlgorithmsMatrix2008}, 
optimization algorithms on manifolds could be used \cite{absilOptimizationAlgorithmsMatrix2008}. 
One may go even further and restrict the optimization to the Grassmann manifold of $m$-dimensional subspaces in $\bbR^K$, denoted $Grass(m,K)$,
see for example \cite[Section 3.4.4]{absilOptimizationAlgorithmsMatrix2008}.
This was for example considered in \cite{liSharpDetectionLowdimensional2025} for linear dimension reduction in Bayesian inverse problems.

There are two main issues with this class of feature maps.
The first is theoretical and lies in the constant $C(\bfX | \calG_m)$, since there is currently no way (even heuristic) to control it.
Worse, if the level set $g^{-1}(\{\bfz\})$ is not \textit{smoothly pathwise-connected}, then we know that $C(\bfX | g(\bfX)= \bfz) = +\infty$, which breaks down the theoretical guarantees of the method.
This can be the case for example when using as feature maps polynomials of degree 2, or when the support of $\mu$ is not convex.
The second issue is that although the coordinates of $g \in \calG_m$ are linear combinations of fixed features $(\Phi_i)_{1\leq j\leq K}$, minimizing $G \mapsto \calJ(G^T \Phi)$ is much more difficult  than in the linear case. 
It does not boil down to solving an eigenvalue problem on some matrix, which is mainly due to the fact that the projection matrix $\Pi_{\nabla g(\bfx)}$ depends on $\bfx$.
Still, despite these two issues, the numerical results in \cite{bigoniNonlinearDimensionReduction2022} give strong positive incentives for using this method, even when connectivity of level sets is not satisfied.

\subsubsection{Feature maps based on diffeomorphisms}

In \cite{zhangLearningNonlinearLevel2019,verdiereDiffeomorphismbasedFeatureLearning2025}, the authors considered for $\calG_m$ nonlinear feature maps defined as the first $m$ components of a class $\calD$ of $\calC^1$-diffeomorphisms of $\bbR^d$,
\begin{equation}
\label{equ:Gm diffeo}
  \calG_m = \{g:\bfx \mapsto (\varphi_1(\bfx), \cdots, \varphi_m(\bfx)) : \varphi \in \calD \},
\end{equation}
where for any $\varphi \in \calD$, it holds  $\varphi,\varphi^{-1} \in \calC^1(\bbR^d, \bbR^d)$.
As with \eqref{equ:Gm vector space}, the variety of possible feature maps $g$ has the advantage of being potentially able to approximate functions with fewer features than \eqref{equ:Gm linear features}.

The main advantage however is that it possibly allows a better control on the Poincar\'e constants.
In particular, it ensures connectivity of the level sets of the feature maps in $\calG_m$, which circumvents the main theoretical issues of \eqref{equ:Gm vector space} and  allows us to derive Poincar\'e inequalities for the conditional random variables $\bfX|g(\bfX)= \bfz$, although the associated Poincar\'e constants are not available in general. 
Moreover, it allows to leverage Poincar\'e inequalities in the feature space $\varphi(\calX)$.
It consists in applying \Cref{prop:poincare based bound} on $u \circ \varphi^{-1} : \varphi(\calX) \rightarrow \bbR$ with features $\varphi_1(\bfX), \cdots, \varphi_m(\bfX)$ which are linear in $\varphi(\bfX)$.
As a result, if one ensures that $\varphi(\bfX) \sim \calN(0, I_d)$, in which case $\varphi$ is called a \textit{normalizing flow} for $\bfX$, then we would obtain $C(\varphi(\bfX) | g(\bfX) = \bfz) = 1$. 
To this end, the authors in \cite{verdiereDiffeomorphismbasedFeatureLearning2025} introduce a penalized optimization problem with a penalty based on some divergence (e.g., the Kullback-Leibler divergence) between  $\varphi(\bfX)$ and $\calN(0, I_d)$, which is a reasonable heuristic, although it gives no theoretical control over the Poincar\'e constant.

The main issue with this class of feature maps is that the set $\calG_m$ has a much more complex structure than the vector space structure in \eqref{equ:Gm vector space}.
It also strongly depends on the choice of parameterization of $\calD$.
As a result, not only the objective function $\calJ$ is highly non-convex, but so is the set on which it is optimized.

We end this section by pointing out that in \cite{verdiereDiffeomorphismbasedFeatureLearning2025} the authors proposed a dimension augmentation strategy, consisting in approximating $u(\bfx)$ by $\bbE(f(g(\bfx,\bfXi))),$ where  $\bfXi \in \bbR^m$    is some random vector, independent on $\bfX$, with connected support. Then the above approach  can be used to construct a map $g : \bbR^{d + m} \to \bbR^{m} $ from  the first $m$ coordinates of a diffeomorphism.  
This approach has two main advantages.
The first is that it enlarges the set $\calD$, thus increasing the expressivity of the admissible feature maps, while preserving the theoretical guarantees.
The second is that it allows to easily define diffeomorphisms  $\varphi: (\bfx, \bfxi) \mapsto (h(\bfx) + \bfxi ,  \bfx) $ on $  \bbR^{d+m}$ from arbitrary functions $h\in\calC^1(\bbR^d, \bbR^m)$, thus  circumventing the theoretical issue of \eqref{equ:Gm vector space}.
However, it also brings some new issues, especially the fact that the approximation of $u$  involves an expectation over $\bfxi$, see the discussion in \cite{verdiereDiffeomorphismbasedFeatureLearning2025}.

\section{Large and small deviations inequalities}
\label{sec:Large and small deviations inequalities}

In the next \Cref{sec:one feature}, we will see in \eqref{equ:norm projection applied 1} that for any scalar valued feature map $g\in\calG_1$ such that $\|\nabla g(\bfX)\|_2>0$ almost surely, we can write $\calJ(g)$ as
\[
  \calJ(g)
  = \Expe{
    \|\nabla g(\bfX)\|_2^{-2}
    \|\nabla u(\bfX)\|_2^2
    \|\Piperp_{\nabla u(\bfX)} \nabla g(\bfX)\|_2^2
  }.
\]
Observing that the term within the above expectation differs only by a factor $\|\nabla g(\bfX)\|_2^2$ from the one in our surrogate $\calL_1(g)$ defined in \eqref{equ:def of J} will be the starting point of our new approach in \Cref{sec:one feature}.
Then, we will obtain theoretical guarantees by leveraging so-called small and large deviations inequalities for the random variable $\| \nabla g({\bf X}) \|_2^2$.
In other words, we require upper bound on both 
\[
  \Proba{\| \nabla g({\bf X}) \|_2^2 \leq t^{-1}}
  \quad \text{and} \quad
  \Proba{\| \nabla g({\bf X}) \|_2^2 \geq t}
\]
as $t \rightarrow +\infty$.
To that end, we summarize in this section some results from  \cite{bobkovSharpDilationtypeInequalities2008,fradeliziConcentrationInequalities$s$concave2009}.
Although these works provided independently similar results, we will be interested in some specific cases studied in \cite{fradeliziConcentrationInequalities$s$concave2009}, and for this reason we will mostly rely on the latter reference.

The results of interest in our framework are sharp large and small deviation inequalities for $s$-concave probability measures with $s\in[-\infty, 1/d]$.
This class of probability measures has been well studied. We refer to the introductions of \cite{borellConvexSetFunctions1975,fradeliziConcentrationInequalities$s$concave2009} and, \cite[Section 3]{bobkovSharpDilationtypeInequalities2008} for a short summary and basic properties.
We also refer to \cite{borellConvexSetFunctions1975,borellConvexMeasuresLocally1974} for a deeper study.
For the sake of simplicity, here we only consider measures admitting a density with respect to the Lebesgue measure, and we introduce them in \Cref{def:s-concave measures} via the characterization recalled in \cite[Introduction]{fradeliziConcentrationInequalities$s$concave2009}.

\begin{definition}[$s$-concave probability measure]
\label{def:s-concave measures}
  Let $\mu$ a probability measure on $\bbR^d$ such that $d\mu(\bfx) = \rho(\bfx) d\bfx$.
  For $s\in[-\infty, 1/d]$, $\mu$ is $s$-concave if and only if $\rho$ is supported on a convex set and is
  $\kappa$-concave with $\kappa=s/(1-sd) \in [-1/d, +\infty]$, meaning
  \begin{equation}
    \rho(\lambda \bfx + (1-\lambda) \bfy)
    \geq (\lambda \rho(\bfx)^{\kappa} + (1-\lambda)\rho(\bfy)^{\kappa})^{1/\kappa}
  \end{equation}
  for all $\bfx,\bfy\in\bbR^d$ such that $\rho(\bfx)\rho(\bfy)>0$ and all $\lambda \in [0,1]$.
  The cases $s\in\{-\infty, 0, 1/d\}$ are interpreted by continuity.
\end{definition}

Let us recall some basic properties of $s$-concave probability measures $\mu$ on $\bbR^d$.
Firstly there exists no such measure for $s>1/d$, and if $\mu$ is $s$-concave, it is also $t$-concave for all $t\leq s$.
Secondly for $s \in (0,1/d]$ the support of $\mu$ is compact.
Also, a measure is $s$-concave with $s=1/d$ if and only if the measure is uniform.
Finally, $s=0$ if and only if $\mu$ is log-concave.
We can also note that $s=-1$ for the Cauchy distribution. 

\begin{remark}
A legitimate question is whether measures described in \Cref{def:s-concave measures} also satisfy Poincar\'e inequalities, which seems to be an interesting property in our context, although we recall that \Cref{prop:poincare based bound} requires Poincar\'e inequalities for the conditional measures $\mu_{\bfX | g(\bfX)=\bfz}$ to be satisfied for all $\bfz \in g(\calX)$.
This question has been answered in 
\cite[Theorem 1.2]{bobkovWeightedIsoperimetricPoincaretype2009}, which states that if $\mu$ is $s$-concave with $-\infty < s \leq 0$ then $\mu$ satisfies a \emph{weighted} Poincar\'e inequality.
More precisely, there exists $C_s > 0$, which depends continuously on $s$, such that for all locally Lipschitz function $h:\bbR^d \rightarrow \bbR$ it holds
\[
  \Expe{(h(\bfX) - \Expe{h(\bfX)})^2}
  \leq C_s \Expe{
    \|\nabla h(\bfX)\|_2^2
    \big(
    \alpha_0
    + s^2 \|\bfX\|^2_2  
    \big)
  },
\]
where $\alpha_0:=\exp(\Expe{\log \|\bfX\|_2}) < +\infty$ is a geometric mean.
In the case $s \in (0, 1/d]$, using the fact that if $\mu$ is $s$-concave then $\mu$ is $t$-concave for all $t\leq s$, we can apply the above inequality with $s=0$ to obtain the following classical Poincaré inequality,
\[
  \Expe{(h(\bfX) - \Expe{h(\bfX)})^2}
  \leq C_0 \alpha_0
  \Expe{
    \|\nabla h(\bfX)\|_2^2
  },
\]
for all locally Lipschitz function $h:\bbR^d \mapsto \bbR$.
In the case $s < 0$, which includes measures with heavy tails such as the generalized Cauchy distribution, only the weighted Poincar\'e inequality holds (with possible refinements for more restricted classes of measures).
Still, it seems possible to leverage it for $s$-concave measures with $s<0$ by adapting the definition of $\calJ$ in \eqref{equ:def of J} and our surrogates from \Cref{sec:one feature,sec:multiple features}.
We leave this to further investigation.
Note that for heavy-tailed measures one should pay more attention to the finiteness of the right-hand side of the weighted Poincar\'e inequality.
\end{remark}

For the class of probability measures described by \Cref{def:s-concave measures}, we will state small and large deviations inequalities for \emph{functions with bounded Chebyshev degree}, using the terminology (which does not seem classical) introduced in \cite{fradeliziConcentrationInequalities$s$concave2009}.
This corresponds to measurable functions $h: \bbR^d \rightarrow \bbR$ such that there exist $ k_h < +\infty$ and $A_h < +\infty$ such that for any straight bounded segment $J\subset \bbR^d$ and any measurable set $I\subset J$,
\begin{equation}
\label{equ:remez inequality multivariate}
  \sup_{\bfx\in J} |h(\bfx)|
  \leq 
  \big(\frac{A_h|J|}{|I|} \big)^{k_h} \sup_{\bfx\in I} |h(\bfx)|.
\end{equation}
Note that such a function $h$ also satisfies this inequality for any $k\geq k_h$ and $A \geq A_h$. 
On the other hand, \eqref{equ:remez inequality multivariate} can be seen as a type of \textit{Remez inequality}, which is more classical, therefore we will stick with the latter phrasing.
Typical functions satisfying {\eqref{equ:remez inequality multivariate}} are polynomials.
Indeed, as stated in \cite{fradeliziConcentrationInequalities$s$concave2009}, if $h=P$ is a multivariate polynomial then $P$ satisfies \eqref{equ:remez inequality multivariate} with
\begin{equation}
\label{equ:k A polynomial}
  k_P = \text{deg}(P), \quad A_P = 4,
\end{equation}
where $\text{deg}(P)$ denotes the total degree of $P$.
A direct consequence is that if $\calG_1$ is a set of polynomial functions with total degree at most $\ell + 1$, then for any $g\in\calG_1$ we have that $h:\bfx \mapsto \|\nabla g(\bfx)\|_2^2$ is a multivariate, real-valued polynomial of degree $2\ell$, which thus satisfies \eqref{equ:remez inequality multivariate} with constants $k_h = 2\ell$ and $A_h = 4$.
This is stated in the following \Cref{prop:k A polynomial features}.
Note that \cite{fradeliziConcentrationInequalities$s$concave2009} provides  other examples of functions satisfying {\eqref{equ:remez inequality multivariate}}, such as exponential polynomials of the form $h : \bfx \mapsto \sum_{j=1}^{N} c_j \exp(i \bfw_j^T \bfx)$ with $c_j \in \bbC$ and $\bfw_j \in \bbR^d$ for $1\leq j\leq N$ and $N\in\bbN^*$, for which it has been shown in \cite{nazarovGeometricKannanLovaszSimonovitsLemma2003} that $h$ satisfies {\eqref{equ:remez inequality multivariate}} with $A_h \leq 316$ and $k_h = N$.

\begin{proposition}
\label{prop:k A polynomial features}
  If $\calG_1$ is a set of polynomials with total degree at most $\ell+1$, then for all $g\in\calG_1$ the function $h : \bfx \mapsto \|\nabla g(\bfx)\|_2^2$ satisfies \eqref{equ:remez inequality multivariate} with constants
  \begin{equation}
  \label{equ:k A polynomial features}
    k = 2\ell,
    \quad A = 4.
  \end{equation}
\end{proposition}

Now, the deviations inequalities in this framework, relying on the analysis of dilatation of the sublevel sets of $\bfx\mapsto |h(\bfx)|$ provided in \cite{fradeliziConcentrationInequalities$s$concave2009,bobkovSharpDilationtypeInequalities2008}, allows us to link the measures of two sublevel sets of the function.
Then, we can choose for example a sublevel set defined by the $\omega$-quantile of $|h(\bfX)|$, denoted by $Q_{|h(\bfX)|}(\omega)$.
We recall that for any real random variable $Z$ the \textit{quantile function} is defined by
\begin{equation}
\label{equ:def quantile function}
  Q_{Z}(\omega):=
  \inf \{z\in\bbR: \Proba{Z \leq z} \geq \omega \},
  \quad \forall \omega\in(0,1).
\end{equation}
In particular, we will only use $q_h := Q_{|h(\bfX)|}(1/2)$, the median of $|h(\bfX)|$.
We can now state in \Cref{prop:small deviation remez uniform,prop:large deviation remez uniform} respectively small and large deviations results.
Note that these results are essentially the ones from \cite{fradeliziConcentrationInequalities$s$concave2009}.
The only difference is that we slightly weakened the bounds in order  to make them easier to leverage in our analysis in further sections, while not changing the rates of decay.

\begin{proposition}[Small deviations]
\label{prop:small deviation remez uniform}
  Let $\bfX$ an absolutely continuous random variable on $\bbR^d$ whose law $\mu$ is $s$-concave with $s\in [-\infty, 1/d]$.
  Let $h: \bbR^d \rightarrow \bbR$ satisfying \eqref{equ:remez inequality multivariate} with constants $k, A < +\infty$.
  Then for all $\varepsilon >0$ it holds
  \begin{equation}
  \label{equ:small deviation remez}
    \Proba{|h(\bfX)| < q_h \varepsilon}
    \leq \underline{\eta}_{A,s} \varepsilon^{1/k}
  \end{equation}
  with $q_h$ defined after \eqref{equ:def quantile function} and $\underline{\eta}_{A,s} := \max(1,(1-2^{-s})s^{-1} A)$.
\end{proposition}
\begin{proof}
  Let $h$ satisfying \eqref{equ:remez inequality multivariate} with constants $k,A$.
  If $h=0$ then the desired inequality holds as $q_h=0$ and $\Proba{|h(\bfX)| < 0} = 0$.
  Otherwise, using the fact that $q_{|h|^{1/k}}^k=q_h$, we have from \cite{fradeliziConcentrationInequalities$s$concave2009} that for all $0<\varepsilon<1$,
  \[
    \Proba{|h(\bfX)| < q_h \varepsilon}
    \leq \Proba{|h(\bfX)| \leq q_h \varepsilon}
    \leq A (1-2^{-s})s^{-1} \varepsilon^{1/k}
    \leq \underline{\eta}_{A,s} \varepsilon^{1/k}.
  \] 
  Then since $\underline{\eta}_{A,s} \varepsilon^{1/k} \geq 1$ whenever $\varepsilon \geq 1$, we have that \eqref{equ:small deviation remez} holds for all $\varepsilon<0$.
\end{proof}

\begin{proposition}[Large deviations]
\label{prop:large deviation remez uniform}
  Let $\bfX$ an absolutely continuous random variable on $\bbR^d$ whose law $\mu$ is $s$-concave with $s\in [-\infty, 1/d]$.
  Let $h: \bbR^d \rightarrow \bbR$ satisfying \eqref{equ:remez inequality multivariate} with constants $k, A < +\infty$.
  Then for all $t>0$ it holds
  \begin{equation}
  \Proba{|h(\bfX)| > q_h t} \leq
  \left\{
  \begin{aligned}
    & \big( 1 - \frac{t^{1/k}-1}{\overline{\eta}_{A,s}} \big)_+^{1/s},&
    & s \in (0,1/d],&
    \\
    & \exp(-\frac{t^{1/k}-1}{\overline{\eta}_{A,s}}),&
    & s=0,&
    \\
    & \overline{\eta}_{A,s} t^{1/s k},&
    & s \in [-\infty, 0),&
  \end{aligned}
  \right.
  \end{equation}
  with $q_h$ defined after \eqref{equ:def quantile function} and
  \begin{equation}
  \overline{\eta}_{A,s} := 
  \left\{
  \begin{aligned}
    &A / (1 - 2^{-s}), &
    & s \in (0,1/d],&
    \\
    & A / \log(2),&
    & s=0,&
    \\
    & \max(1,A / (2^{-s} - 1))^{-1/s}, &
    & s \in [-\infty, 0).&
  \end{aligned}
  \right.
  \end{equation}
\end{proposition}
\begin{proof}
  See \Cref{subsec:proof prop:large deviation remez uniform}
\end{proof}

Although the analysis in the next sections can be performed using medians, it will be more convenient to rely on $L^p_{\mu}$ norms instead, especially with $p\geq 1$.
It turns out that for functions satisfying \eqref{equ:remez inequality multivariate}, the median can be explicitly lower and upper bounded using $L^p_{\mu}$ norms.
This is stated in \Cref{prop:bounded uniform quantiles}.
Note that the upper bound on the median comes from the Markov's inequality, and that the lower bound is only a restatement of results from \cite{fradeliziConcentrationInequalities$s$concave2009} that directly leverage large deviation inequalities.
Note also that for $s\in(-1/k,1/d)$, these results imply that if $h$ satisfies \eqref{equ:remez inequality multivariate} then $h\in L^p_{\mu}$ for all $1\leq p<-1/sk$.

\begin{proposition}
\label{prop:bounded uniform quantiles}
  Let $\bfX$ an absolutely continuous random variable on $\bbR^d$ whose law $\mu$ is $s$-concave with $s\in(-1/k, 1/d]$.
  Let $h:\bbR^d \rightarrow \bbR$ not identically zero satisfying \eqref{equ:remez inequality multivariate} with constants $k$ and $A$.
  Then, for any $p\geq 1$ such that $-1 < spk$ it holds 
  \begin{equation}
    2^{\frac{1}{p}}
    \geq
    \frac{q_h}{\|h\|_{L^p_{\mu}}} 
    \geq
    A^{-k}
    \left\{
    \begin{aligned}
      &(1-2^{-s})^k, 
      & s\in(0,1/d],\\
      & (3pk)^{-k}, 
      & s \in [0,1/d],\\
      &\Big(1 - \frac{(2^{-s} - 1)^{\frac{1}{s}}}{1 + (spk)^{-1}}\Big)^{-\frac{1}{p}},
      & s\in(-1/k, 0). \\
    \end{aligned}
    \right.
  \end{equation}
\end{proposition}
\begin{proof} 
  See \Cref{subsec:proof prop:bounded uniform quantiles}.
\end{proof}

Note that in the proof of \Cref{prop:bounded uniform quantiles} we provide bounds that holds also for $p\in(0,1)$, however in this work we will mainly consider results for $p\geq 1$.
Let us now consider two examples for a $s$-concave measure $\mu$ and a polynomial $h:\bbR^d \mapsto \bbR$ of total degree $k$, so that $h$ satisfies \eqref{equ:remez inequality multivariate} with constants $k$ and $A=4$.
Firstly, let $s=1/d$, meaning that $\mu$ is uniform on a convex set $\calX \in \bbR^d$, then the constants in \Cref{prop:bounded uniform quantiles} satisfy for all $p\in[1,+\infty]$,
\[
  \underline{\eta}_{4,\frac{1}{d}} \leq 4 \log(2),
  \quad
  \overline{\eta}_{4,\frac{1}{d}} 
  \leq 8d,
  \quad 
  \frac{1}{(8d)^k} \|h\|_{L^p_{\mu}} 
  \leq q_h
  \leq  2^{1/p}\|h\|_{L^p_{\mu}}.
\]
Note that the above inequalities are almost independent on $p$ but strongly depend on $d$.
Such bound for $s>0$ will be especially useful in \Cref{prop:small deviation sigma m} in which we leverage control in $L^{\infty}_{\mu}$.
Secondly, and more generally, let $s\geq 0$, meaning that $\mu$ is log-concave.
Then for any $p \in [1,+\infty)$ it holds 
\[
  \underline{\eta}_{4,0} = 4 \log(2),
  \quad 
  \overline{\eta}_{4,0} = \frac{4}{\log(2)},
  \quad
  \frac{1}{(3pk)^k} \|h\|_{L^{p}_{\mu}} 
  \leq  q_h 
  \leq 2^{1/p} \|h\|_{L^{p}_{\mu}} .
\]
The main difference with the case $s=1/d$ is that the above bounds do not depend on $d$ or $s$ but strongly depend on $p$.
This dependency in $p$ comes from the fact that $0$-concave measures may have unbounded support, as for the Gaussian measure.
It is important to note that for $s>0$, depending on the value $p$ for which we have control $L^{p}_{\mu}$ norm, one may prefer to use one bound or the other.
More precisely for $s>0$, if $3pk (1-2^ {-s}) < 1$, then one would prefer to use the general bound that does not depend on $s$.

Let us end this section with a last small deviation result that we will leverage in \Cref{sec:multiple features}.
More precisely, we state in \Cref{prop:small deviation sigma m} a small deviation result on singular values of matrix-valued polynomial functions.
This result mainly relies on applying \Cref{prop:small deviation remez uniform} to functions of the form $\bfx \mapsto \det(M(\bfx))$ for some matrix-valued polynomial $M:\bbR^d \mapsto \bbR^{m\times m}$, so that $\bfx \mapsto \det(M(\bfx))$ is also polynomial.

\begin{proposition}
\label{prop:small deviation sigma m}
  Assume that $\bfX$ is an absolutely continuous random variable on $\bbR^d$ whose law $\mu$ is $s$-concave with $s\in(0,1/d]$.
  Let $M:\bbR^d \rightarrow \bbR^{m\times m}$ a polynomial with total degree at most $\ell$.
  Then for all $\varepsilon>0$ it holds
  \begin{equation}
  \label{equ:small deviation sigma m}
    \Proba{\sigma_m(M(\bfX)) < q_{\det(M)} \varepsilon}
    \leq \underline{\eta}_{M, s} \varepsilon^{1/\ell m},
  \end{equation}
  with $q$ defined after \eqref{equ:def quantile function} and
  \[
    \underline{\eta}_{M, s}
    := 
    \underline{\eta}_{4,s} 
    (2m^{-1})^{\frac{1}{2\ell}}
    \sup_{\bfx\in\calX} \|M(\bfx)\|_F^{\frac{m-1}{\ell m}}.
  \]
\end{proposition}
\begin{proof}
  Let $M : \bbR^d \mapsto \bbR^{m\times m}$ a polynomial with total degree at most $\ell$.
  If $M=0$ then the desired result holds as $q_{\det(M)}=0$ and $\Proba{\sigma_m(M(\bfX))<0}=0$.
  Otherwise, from \cite{piazzaUpperBoundCondition2002,gungorErratumUpperBound2010} it holds $\sigma_m(M(\bfX)) \geq \kappa_{M}^{-1}|\det(M(\bfX))|$ where
  $
    \kappa_{M} := 
    (m-1)^{-\frac{m-1}{2}}
    \sup_{\bfx\in\calX}
    \|M(\bfx)\|_F^{m-1}.
  $
  Denoting $q=q_{\det(M)}$, the previous inequality then implies for all $\varepsilon>0$,
  \[
    \Proba{\sigma(M(\bfX)) < q \varepsilon} 
    \leq \Proba{|\det(M(\bfX))| \leq q \kappa_{M}\varepsilon}.
  \]
  Now, since $\det: \bbR^{m\times m} \rightarrow \bbR$ is polynomial with total degree at most $m$ and $\bfx \mapsto M(\bfx)$ is polynomial with total degree at most $\ell$, $\bfx \mapsto \det(M(\bfx))$ is polynomial with total degree at most $\ell m$, thus it satisfies \eqref{equ:remez inequality multivariate} with constants $k=\ell m$ and $A=4$.
  Hence, from \Cref{prop:small deviation remez uniform} we obtain for all $\varepsilon>0$,
  \[
    \Proba{\sigma(M(\bfX)) < q \varepsilon}
    \leq 
    \Proba{|\det(M(\bfX))| < q \kappa_{M} \varepsilon}
    \leq \underline{\eta}_{4,s} \kappa_M^{1/\ell m} \varepsilon^{1/\ell m}.
  \]
  Finally, studying the variations of $x \mapsto x(x-1)^{-1 + 1/x}$ we obtain that $m^{-1} \leq (m-1)^{-(m-1)/m}\leq 2 m^{-1}$, hence $\kappa_M^{1/\ell m} \leq (2m^{-1})^{\frac{1}{2\ell}}\sup_{\bfx\in\calX} \|M(\bfx)\|_F^{\frac{m-1}{\ell m}}$, which yields the desired result.
\end{proof}

Let us emphasize four points concerning \Cref{prop:small deviation sigma m}.
Firstly, the constant $\underline{\eta}_{M,s}$ depends on $M$, while $\underline{\eta}_{A,s}$ from \Cref{prop:small deviation remez uniform} only depends on $A$ and $s$.
Secondly, we only state the result for polynomials, because if $\bfx \mapsto M(\bfx)$ is polynomial then $\bfx \mapsto \det(M(\bfx))$ is also polynomial, which satisfies \eqref{equ:remez inequality multivariate}.
However, we were not able to identify a more general class of functions such that $\bfx \mapsto \det(M(\bfx))$ satisfies \eqref{equ:remez inequality multivariate}.
Thirdly, we restricted our analysis to the case $s\in(0,1]$ for the sake of simplicity.
One could generalize to the case $s\leq 0$ by additionally considering large deviation inequalities on $\|M(\bfX)\|_F^2$.
However, this would result in a more complex statement, whose complexity would propagate to \Cref{prop:suboptimality uniform m}.
Fourthly, we can expect the exponent $1/\ell m$ in \eqref{equ:small deviation sigma m} to be sharp.
For example one can consider $X \sim \calU((0,1))$ and $M : \bbR \mapsto \bbR^{2\times 2}$ with $M(x)_{1,1}=M(x)_{2,2} = x$, $M(x)_{1,0}=0$ and $M(x)_{0,1} = 1$ such that $m=2$ and $\ell=1$.
One can then show that $\sigma_2(M(x)) \sim x^2$ as $x\rightarrow 0$.
As a result, $\Proba{\sigma_2(M(X)) \leq \varepsilon} \gtrsim \varepsilon^{1/2}$ as $\varepsilon \rightarrow 0$.

\section{The case of a single feature}
\label{sec:one feature}

In this section we study the case $m=1$ where we search for a scalar-valued feature map $g$.
We  introduce a new quantity, for any $g\in\calG_1$,
\begin{equation}
\label{equ:def L one feature}
  \calL_1(g)
  := \Expe{\|\nabla u(\bfX)\|_2^2 \, \|\Piperp_{\nabla u(\bfX)} \nabla g(\bfX) \|_2^2},
\end{equation}
whose minimizer, if it exists, shall be a candidate feature map.
The idea behind considering this very surrogate lies in the following \Cref{lem:norm projection}.

\begin{lemma}
\label{lem:norm projection}
  For all $\bfv, \bfw \in \bbR^d$, it holds
  \begin{equation}
    \|\bfw\|^2_2 \|\Pi_{\bfw^{\perp}} \bfv \|^2_2 = \|\bfv\|^2_2 \|\Pi_{\bfv^{\perp}} \bfw\|^2_2.
  \end{equation}
\end{lemma}
\begin{proof}
  It holds
  \[
    \|\bfw\|^2_2 \|\Piperp_{\bfw}\bfv \|^2_2
    = \|\bfw\|^2_2\|\bfv\|^2_2 - \|\bfw\|^2_2\innerp{\frac{\bfw}{\|\bfw\|_2}}{\bfv}^2_2
    = \|\bfw\|^2_2\|\bfv\|^2_2 - \innerp{\bfw}{\bfv}^2_2
    = \|\bfv\|^2_2\|\Piperp_{\bfv}\bfw \|^2_2.
  \]
\end{proof}

Thus, if $\calG_1$ satisfies \Cref{assump:full rank jacobian as}, then for any $g\in\calG_1$, taking $\bfw=\nabla g(\bfX)$ and $\bfv=\nabla u(\bfX)$ in the above lemma yields the alternative expressions for $\calJ(g)$ and $\calL_1(g)$,
\begin{equation}
\label{equ:norm projection applied 1}
\begin{aligned}
  \calJ(g)
  &= \Expe{
    \|\nabla g(\bfX)\|_2^{-2}
    \|\nabla u(\bfX)\|_2^2
    \|\Piperp_{\nabla u(\bfX)} \nabla g(\bfX)\|_2^2
    }, \\
  \calL_1(g) 
  &= \Expe{
    \|\nabla g(\bfX)\|_2^2
    \|\Piperp_{\nabla g(\bfX)} \nabla u(\bfX)\|_2^2
    }.
\end{aligned}
\end{equation}

In the first and second right-hand sides we recognize the terms within the expectation in respectively $\calL_1(g)$ and $\calJ(g)$, up to a factor $\|\nabla g(\bfX)\|_2^2$.
Thus, controlling elements of 
\begin{equation}
\label{equ:def set norm grad features}
  \calK_1 := \{\bfx \mapsto \|\nabla g(\bfx)\|_2^2 \,: \, g\in\calG_1\}
\end{equation}
is the key point for controlling the suboptimality of $\calL_1$.
For example for linear features $g(\bfx) = G^T \bfx$ with $G^T G = 1$ we have that $\|\nabla g(\bfX)\|_2 = 1$ almost surely thus $\calL_1(g) = \calJ(g)$.
A first element of control is given by \Cref{assump:full rank jacobian as}, which implies that $\|\nabla g(\bfX)\|_2^2>0$ almost surely in the case $m=1$.
Ideally, we would like to uniformly lower and upper bound elements of $\calK_1$ in order to obtain a strong  link between $\calJ(g)$ and $\calL_1(g)$.
However, when this is not possible, for example when considering nonlinear polynomials in $\calG_1$, we can only rely on a weaker control.
In our case, we will rely on concentration-type inequalities detailed in \Cref{sec:Large and small deviations inequalities}, which can be leveraged assuming that $\calG_1$ satisfies the following.

\begin{assumption}
\label{assump:bounded uniform A k}
  For all $g\in\calG_1$, the function $\bfx \mapsto \|\nabla g(\bfx)\|_2^2$ satisfies \eqref{equ:remez inequality multivariate} with constants $k<+\infty$ and $A<+\infty$.
\end{assumption}
One shall keep in mind that, as stated in \Cref{prop:k A polynomial features} and the discussion after, if $\calG_1$ contains only polynomial feature maps of total degree at most $\ell+1$, then $\calG_1$ satisfies the above \Cref{assump:bounded uniform A k} with $k = 2\ell$ and $A=4$.
Additionally, if $\calG_1$ does not contain any constant polynomial, then it also satisfies \Cref{assump:full rank jacobian as}. 
It is also worth noting that feature maps based on trigonometric functions as considered in \cite{romorKernelbasedActiveSubspaces2022} also satisfy \Cref{assump:bounded uniform A k}, as stated in \Cref{rem:trigonometric features}. 
However, the main problem is that the associated constant $k$ can be expected to be large.

\begin{remark}[Trigonometric features]
\label{rem:trigonometric features}
  In \cite{romorKernelbasedActiveSubspaces2022} the authors considered a set of feature maps $\calG_m$ as in \eqref{equ:Gm vector space sphere} with $\Phi_j(\bfx) = \cos(\bfw_j^T \bfx + b_j)$ for some $\bfw_j\in\bbR^d$ and $b_j\in\bbR$ for all $1\leq j\leq K$.
  If $m=1$ and if $(\bfw_j)_{1\leq j\leq K}$ are orthonormal vectors in $\bbR^d$, then one can show that for $g(\bfx)=G^T \Phi(\bfx)$ with $G\in\bbR^K$ it holds $\| \nabla g(\bfx) \|_2^2 = \sum_{j=-N}^{N} c_j \exp(i2\bfw_j^T \bfx)$ with $c_0=\|G\|_2^2/2$, $\bfw_0=0$, $c_j = \bar c_{-j} = G_j^2 \exp(i2b_j)/4$ and $\bfw_j=-\bfw_{-j}$ for all $1\leq j\leq N$.
  As a result, in view of the discussion in \Cref{sec:Large and small deviations inequalities} on exponential polynomials, $\calG_1$ satisfies \Cref{assump:bounded uniform A k} with constant $k\leq 2K+1$ and $A\leq 316$.
  Now if $(\bfw_j)_{1\leq j\leq K}$ are not orthogonal, then one can still show that $\bfx \mapsto \|\nabla g(\bfx)\|_2^2$ is an exponential polynomial, however the number of terms seems to scale as $K^2$.
  Both cases are to compare with polynomial feature maps with total degree at most $\ell+1$, for which $\calG_1$ satisfies \Cref{assump:bounded uniform A k} with $k=2\ell$ and $A=4$ and for which $K \leq \binom{d}{\ell+1}$.
\end{remark}

The structure of the section is as follows.
First in \Cref{subsec:sub optimality} we provide results on the suboptimality of $\calL_1$.
Then in \Cref{subsec:minimizing the surrogate} we discuss on the minimization of $\calL_1$.

\subsection{Suboptimality}
\label{subsec:sub optimality}

In this section, we discuss on the suboptimality of our surrogate $\calL_1$ defined in \eqref{equ:def L one feature} regarding the Poincaré inequality error bound $\calJ$ defined in \eqref{equ:def of J}.

Let us start by investigating the special case when $u$ can be described with exactly one feature in $\calG_1$, in other words when $\calJ(g) = 0$ for some $g\in\calG_1$.
In this case, we can show that minimizing $\calL_1$ is equivalent to minimizing $\calJ$, as stated in \Cref{prop:exact recovery}.
Note that this result relies on the fact that there is no more than $1$ feature to learn.

\begin{proposition}
\label{prop:exact recovery}
  Assume $\calG_1$ satisfies \Cref{assump:full rank jacobian as}.
  Then for all $g\in\calG_1$ we have $\calJ(g)=0$ if and only if $\calL_1(g)=0$.
\end{proposition}
\begin{proof}
  Let us start by assuming that $\calJ(g) = 0$.
  The definition of $\calJ$ then yields that $\nabla u(\bfX) \in \spanv{\nabla g(\bfX)}$ almost surely.
  If $\nabla u(\bfX) \neq 0$, then since we are in the one feature case, this implies that $\nabla g(\bfX) \in \spanv{\nabla u(\bfX)}$.
  Otherwise, we have $\|\nabla u(\bfX)\|_2=0$.
  Hence in both cases, we have
  \[
    \|\nabla u(\bfX)\|_2^2 
    \|\Piperp_{\nabla u(\bfX)} \nabla g(\bfX)\|_2^2
    =0
  \]
  almost surely.
  As a result $\calL_1(g) = 0$.
  Now for the converse statement, let us assume that $\calL_1(g) = 0$.
  The definition of $\calL_1$ then yields that almost surely, $\nabla  u(\bfX)=0$ or $\nabla g(\bfX) \in \spanv{\nabla u(\bfX)}$.
  Then by \Cref{assump:full rank jacobian as} and since we are in the one feature case, this implies that $\nabla u(\bfX) \in \spanv{\nabla g(\bfX)}$ almost surely, which implies that 
  \[
    \|\Piperp_{\nabla g(\bfX)} \nabla u(\bfX)\|_2^2
    =0,
  \]
  almost surely.
  As a result, $\calJ(g) = 0$.
\end{proof}

Now in more general cases we cannot expect that $\min_{g\in \calG_1}\calJ(g) = 0$, and we would want to quantify the suboptimality of $\calL_1$.
An ideal setting, considered in the following \Cref{prop:sub optimality bi lip}, is when $g$ is bi-Lipschitz  with constants uniformly controlled on $\calG_1$.
Note that this is for example the case when  $\calG_1$ is defined by \eqref{equ:Gm diffeo} with $\calD$ containing only bi-Lipschitz diffeomorphisms whose associated constants are uniformly bounded over $\calD$.

\begin{proposition}
\label{prop:sub optimality bi lip}
  Assume that there exist $c>0$ and $C < +\infty$ such that for all $g\in\calG_1$ and $\bfx\in\calX$ it holds $c \leq \|\nabla g(\bfx)\|_2^2 \leq C$.
  Then for all $g\in\calG_1$,
  \begin{equation}
    c\calJ(g)
    \leq 
    \calL_1(g) 
    \leq
    C \calJ(g).
  \end{equation}
  Moreover, if $g^*$ is a minimizer of $\calL_1$ over $\calG_1$, it yields
  \begin{equation}
    \calJ(g^*) 
    \leq (C / c) \inf_{g\in \calG_1} \calJ(g).
  \end{equation}
\end{proposition}
\begin{proof}
  The result follows directly from \eqref{equ:norm projection applied 1}.
\end{proof}

In the bi-Lipschitz setting, $\calJ$ and $\calL_1$ are equivalent in some sense.
However, this setting is very restrictive, especially regarding the lower bound on the gradient.
For example if $\calG_m$ is taken as \eqref{equ:Gm vector space sphere} with $\Phi:\bbR^d \rightarrow \bbR^K$ and $K>d$, then for fixed $\bfx\in\calX$ we have rank$(\nabla \Phi(\bfx)) \leq d$, thus there exists $G\in\bbR^K \setminus \{0\}$ such that $G^T \nabla \Phi(\bfx) = 0 = \nabla g(\bfx)$ where $g=G^T \Phi$.
Thus in this case, there cannot exist bi-Lipschitz constants uniformly bounded on $\calG_1$. 
For this reason, the theoretical guarantees on $\calL_1$ will heavily rely in general  on concentration-type inequalities in order to upper bound the probabilities of $(\|\nabla g(\bfX)\|_2^2\geq \beta)$ and $(\|\nabla g(\bfX)\|_2^2\leq \beta^{-1})$ as $\beta \rightarrow +\infty$.
The first one is referred to as \textit{large deviation} probability and widely known concentration inequalities (e.g., Chernoff's) provide the desired bounds.
The second one is referred to as \textit{small deviation} probability, which is much less studied than  large deviations.

We now introduce quantities that uniformly control the smallest and largest elements of a set $\calK$ of functions, with respect to the $L^r_{\mu}$ norm for $r\geq 1$,
\begin{equation}
\label{equ:def nu}
  \underline{\nu}_{\calK, r} 
  := \inf_{h\in\calK} \Expe{|h(\bfX)|^r}^{1/r},
  \quad
  \overline{\nu}_{\calK, r}
  := \sup_{h\in\calK} \Expe{|h(\bfX)|^r}^{1/r}.
\end{equation}
Let us first consider the control of $\calJ$ by $\calL_1$ uniformly on $\calG_1$.
We state such a result in \Cref{prop:bound J by L}, and we recall that from \Cref{prop:k A polynomial features}, if $\calG_1$ contains only polynomials of total degree at most $\ell+1$ then $\calG_1$ satisfies \Cref{assump:bounded uniform A k} with $k = 2\ell$ and $A=4$.
Note that we only consider $s$-concave measures with $s>-1/k$ because \Cref{prop:bounded uniform quantiles} then ensures that $\nu_{\calK_1, 1} < +\infty$, therefore $\calL_1(g)<+\infty$.
On the other hand if $s\leq -1/k$ there exists heavy-tailed $s$-concave distributions such that some polynomials of total degree $k$ are unbounded in $L_{\mu}^1$.
Note also that, for the sake of readability, we state our results assuming the function $u$ satisfies $\Expe{\|\nabla u(\bfX)\|_2^{2p_u}} \leq 1$ for some $p_u >1$.
Equivalent results can be obtained without this assumption by considering $\tilde u: \bfx \mapsto \alpha^{-1/2} u(\bfx)$ with $\alpha := \Expe{\|\nabla u(\bfX)\|_2^{2p_u}}^{1/p_u}$, with associated quantities $\widetilde \calJ(g)$ and $\widetilde \calL_1(g)$, and using $\calJ(g) = \alpha \widetilde \calJ(g)$ and $\calL_1(g) = \alpha \widetilde \calL_1(g)$. 

\begin{proposition}
\label{prop:bound J by L}
  Assume $\bfX$ is an absolutely continuous random variable on $\bbR^d$ whose distribution is $s$-concave with $s\in (-1/k, 1/d]$.
  Assume that $\calG_1$ satisfies \Cref{assump:full rank jacobian as,assump:bounded uniform A k} and that $\Expe{\|\nabla u(\bfX)\|_2^{2p_u}} \leq 1 $ for some $p_u > 1$.
  Then for $p^{-1} := 1 - p_u^{-1}$ and for all $p_1\geq 1$ such that $-1 < skp_1 \leq +\infty$, it holds for all $g\in\calG_1$,
  \begin{equation}
  \label{equ:bound J by L}
    \calJ(g) \leq 
    \gamma_1 
    \underline{\nu}_{\calK_1,p_1}^{-\frac{1}{1+pk}}
    \calL_1(g)^{\frac{1}{1+pk}},
  \end{equation}
  with $\underline{\nu}_{\calK_1,p_1}$ defined in \eqref{equ:def nu}, with $\calK_1$ defined in \eqref{equ:def set norm grad features} and with $\gamma_1$ defined by
  \[
    \gamma_1 := 
    \left\{
    \begin{aligned}
      & 2 \big(
        \underline{\eta}_{A,s} 
        A \min\{3kp_1,(1-2^{-s})^{-1}\}
        \big)^{\frac{k}{1+pk}}
      \quad & s\in(0,1/d], \\
      & 2 \big(
        3\underline{\eta}_{A,0} 
        A k p_1
        \big)^{\frac{k}{1+pk}}
      \quad & s=0, \\
      & 2 \big(
        \underline{\eta}_{A,s}
        A (1- \frac{(2^{-s}-1)^{1/s}}{1+ 1/skp_1} )^{\frac{1}{kp_1}}
        \big)^{\frac{k}{1+pk}}
      \quad & s\in(-1/k, 0), \\
    \end{aligned}
    \right.
  \]
  which depends only on $A,s,k,p_u$ and $p_1$, with $\underline{\eta}_{A,s}$ from \Cref{prop:small deviation remez uniform}.
\end{proposition}
\begin{proof}
See \Cref{subsec:proof prop:bound J by L}.
\end{proof}

We can first note that if $s>0$ then $\bfX$ has compact support, thus we can apply \Cref{prop:bound J by L} with $p_u = +\infty$ and $p=1$.
Also, it is important to note that the bound \eqref{equ:bound J by L} depends on the parameter $s$ of the distribution of $\bfX$ only through the constant $\gamma_1$.
On the other hand, it heavily depends on the constant $k$.
Now we consider the control of $\calL_1$ by $\calJ$ uniformly on $\calG_1$ in \Cref{prop:bound L by J}.

\begin{proposition}
\label{prop:bound L by J}
  Assume $\bfX$ is an absolutely continuous random variable on $\bbR^d$ whose distribution is $s$-concave with $s\in [-1/k, 1/d]$.
  Assume that $\calG_1$ satisfies \Cref{assump:full rank jacobian as,assump:bounded uniform A k} and that $\Expe{\|\nabla u(\bfX)\|_2^{2p_u}} \leq 1 $ for some $p_u > 1$ such that $1-p_u^{-1} > -sk$.
  Then for all $p_1\geq 1$ such that $1-p_u^{-1} > p_1^{-1} > -sk$, for $r^{-1} := 1 - p_u^{-1} - p_1^{-1}>0$ and for all $g\in\calG_1$,
  \begin{equation}
  \label{equ:bound L by J}
    \calL_1(g) 
    \leq \gamma_2\overline{\nu}_{\calK_1,p_1}
    \left\{
    \begin{aligned}
      & \calJ(g), \quad 
      & \gamma_2 =& 2 \overline{\eta}_{A,s}^{k} \quad
      & s\in(0,1/d], 
      \\
      & \calJ(g) (1 + |\log(\calJ(g))|^k ),  \quad 
      & \gamma_2 =& 2 (\overline{\eta}_{A,0} r)^{k} \quad
      & s=0, 
      \\
      & \calJ(g)^{\frac{1}{1-srk}}, \quad
      & \gamma_2 =& 4\overline{\eta}_{A,s}^{1/r} \quad 
      & s < 0, 
      \\
    \end{aligned}
    \right.
  \end{equation}
  with $\overline{\nu}_{\calK_1, p_1}$ defined in \eqref{equ:def nu}, with $\calK_1$ defined in \eqref{equ:def set norm grad features} and where $\gamma_2$ only depends on $A,s,k,p_u$ and $p_1$, with $\overline{\eta}_{A,s}$ from \Cref{prop:large deviation remez uniform}.
\end{proposition}
\begin{proof}
  See \Cref{subsec:proof prop:bound L by J}.
\end{proof}

It is worth noting that the resulting bound strongly depends on the parameter $s$ of the law of $\bfX$, contrarily to the converse statement from \Cref{prop:bound J by L}.
Let us briefly discuss on the different cases for $s$.
First, the case $s\in (0,1/d]$ yields the best bound in terms of $\calJ(g)$, and depends on $k$ only via the constant $\gamma_2$.
This could be expected since in that case $\bfX$ has compact support, thus we can take $p_u=+\infty$ and $r=1$.
Note however that the constant $\gamma_2$ grows as $d^k$ in this case, which we suspect to be sharp when $\calJ(g)$ tends to $0$.
Then, the case $s=0$, which we recall is equivalent to log-concave distributions, yields a slightly worse bound in terms of $\calJ(g)$ due to the logarithmic term raised to the power $k$, which can still be considered as a rather weak dependency on $k$ overall.
Note that the constant $\gamma_2$ does not depend on $d$ in this case, therefore this bound may be sharper than the bound for $s>0$, especially when $d$ is large and $\calJ(g)$ is not too small. 
Finally, the case $s<0$ yields the worst bound in terms of $\calJ(g)$, whose dependency on $k$ is roughly the same as in \eqref{equ:bound J by L}, up to a factor $s$.

Now, we provide in  \Cref{prop:sub optimality} a control on the suboptimality of $\calL_1$.
 
\begin{proposition}
\label{prop:sub optimality}
  Assume $\bfX$ is an absolutely continuous random variable on $\bbR^d$ whose distribution is $s$-concave with $s\in (-1/k, 1/d]$.
  Assume that $\calG_1$ satisfies \Cref{assump:full rank jacobian as,assump:bounded uniform A k} and that $\Expe{\|\nabla u(\bfX)\|_2^{2p_u}} \leq 1 $ for some $p_u > 1$ such that $1-p_u^{-1} > -sk$.
  Then for all $p_1\geq 1$ such that $1-p_u^{-1} > p_1^{-1} > -sk$, for $p^{-1}:= 1-p_u^{-1}$ and $r^{-1} := p^{-1}- p_1^{-1} >0$, any $g^*$ minimizer of $\calL_1$ over $\calG_1$ satisfies
  \begin{equation}
  \label{equ:suboptimality}
    \calJ(g^*) \leq \gamma_3
    \big(
      \frac{\overline{\nu}_{\calK_1,p_1}}{\underline{\nu}_{\calK_1,p_1}}
    \big)^{\frac{1}{1+pk}} 
    \left\{
    \begin{aligned}
      &\inf_{g\in \calG_1} \calJ(g)^{\frac{1}{1+pk}} 
      \quad &s\in(0,1/d], \\
      &\inf_{g\in \calG_1}\calJ(g)^{\frac{1}{1+pk}} 
      (1 +|\log(\calJ(g))|^k)^{\frac{1}{1+pk}} 
      \quad &s=0, \\
      &\inf_{g\in \calG_1} \calJ(g)^{\frac{1}{1+pk}\frac{1}{1-srk}} \quad &s < 0, \\
    \end{aligned}
    \right.
  \end{equation}
  with $\underline{\nu}_{\calK_1, p_1}$ and $\overline{\nu}_{\calK_1, p_1}$ defined in \eqref{equ:def nu}, with $\calK_1$ defined in \eqref{equ:def set norm grad features} and with $\gamma_3 := \gamma_1 \gamma_2^{\frac{1}{1+pk}}$ where $\gamma_1$ and $\gamma_2$ are defined in  \Cref{prop:bound J by L,prop:bound L by J} respectively and depend only on $A,s,k,p_u$ and $p_1$.
\end{proposition}
\begin{proof}
  Let $g^* \in \calG_1$ be a minimizer of $\calL_1$ over $\calG_1$.
  Then for all $g\in\calG_1$ we have by \Cref{prop:bound J by L} that $\calJ(g^*) \leq \gamma_1 (\underline{\nu}_{\calK_1,p_1}^{-1} \calL_1(g^*))^{\frac{1}{1 + pk}} \leq \gamma_1 (\underline{\nu}_{\calK_1,p_1}^{-1}\calL_1(g))^{\frac{1}{1 + pk}}$.
  Now, using \Cref{prop:bound L by J} and taking the infimum over $g\in\calG_1$ yields the desired result.
\end{proof}

Let us now put a stress on the special case $s=1/d$, meaning when $\bfX$ is uniformly distributed on a compact convex set of $\bbR^d$, with polynomial feature maps with bounded total degree such that $\Expe{\|\nabla g(\bfX)\|_2^2}=1$ for all $g\in\calG_1$.
Note that such feature maps can be obtained by taking $\Phi:\bbR^d \rightarrow \bbR^K$ polynomial and $R = \Expe{\nabla \Phi(\bfX)^T \nabla \Phi(\bfX)}$ in \eqref{equ:Gm vector space sphere}.
In this case, our previous results become the much simpler \Cref{cor:suboptimality uniform} stated below.

\begin{corollary} 
\label{cor:suboptimality uniform}
  Assume that $\bfX$ is uniformly distributed on a compact convex subset of $\bbR^d$.
  Assume that every $g\in \calG_1$ is a non-constant polynomial of total degree at most $\ell+1$ such that $\Expe{\|\nabla g(\bfX)\|_2^2}=1$.
  Assume that $\sup_{\bfx\in\calX}\|\nabla u(\bfx)\|_2^2 \leq 1$.
  Then for all $g\in\calG_1$,
  \begin{equation}
    \gamma_2^{-1} \calL_1(g) 
    \leq \calJ(g) 
    \leq \gamma_1 \calL_1(g)^{\frac{1}{1+2\ell}},
    \quad \gamma_1 \leq 2(32\min\{3 \ell, d\})^{\frac{2\ell}{1+2\ell}},
    \quad \gamma_2 \leq 2(8d)^{2\ell}.
  \end{equation}
  Moreover, for $g^*$ a minimizer of $\calL_1$ over $\calG_1$ it holds
  \begin{equation}
  \label{equ:def gamma 3 uniform}
    \calJ(g^*) \leq \gamma_3 
    \inf_{g\in\calG_1} \calJ(g)^{\frac{1}{1+2\ell}},
    \quad 
    \gamma_3 \leq 4(256 d \min\{3\ell, d\})^{\frac{2\ell}{1+2\ell}}.
  \end{equation}
\end{corollary}
\begin{proof}
  See \Cref{subsec:proof cor:suboptimality uniform}.
\end{proof}

Let us end this section by discussing on the sharpness of \eqref{equ:bound J by L}, which is the main cause of discrepancy between $\calJ$ and $\calL_1$, especially due to the exponent $1 / (1+pk)$.
For the sake of simplicity, we focus on the case $p_u=+\infty$.
Since the exponent $1/(1+k)$ comes from the exponent in the right-hand side of \eqref{equ:small deviation remez}, we may ask whether this exponent can be improved.
The answer is no since for $X \sim \calU((0,1))$, $k\in \bbN^*$ and $\varepsilon \in (0,1)$ it holds $\Proba{X^k \leq \varepsilon} = \varepsilon^{1/k}$.
Similarly, the factor $k$ in the exponent in \Cref{prop:bound J by L} cannot be improved in general, as shown in \Cref{prop:sharpness of rate with poly degree}.

\begin{proposition}
\label{prop:sharpness of rate with poly degree}
  Let $\bfX \sim \calU((0,1)^2)$, $\ell\in\bbN^*$ and $u(\bfx) = x_1$.
  Consider the set of feature maps
  \[
    \calG_1 = \{g^a:~ 0 \leq a\leq 1\},
    \quad 
    g^a(\bfx) = \kappa_a^{-1}(x_1^{\ell+1} + a^{\ell} x_2^{\ell+1}),
    \quad 
    \kappa_a := (\ell+1) \Expe{x_1^{2\ell} + a^{2\ell} x_2^{2\ell}},
  \]
  so that $\calG_1$ satisfies \Cref{assump:full rank jacobian as,assump:bounded uniform A k} and $\Expe{\|\nabla g^a(\bfX)\|_2^2} = 1$ for any $0 \leq a\leq 1$.
  Then it holds $\calJ(g^a) \geq a/4$ and $\calL_1(g^a) = a^{2\ell} / (1 + a^{2\ell})$ for $a>0$ small enough.
  Moreover, for $a \rightarrow 0$ it holds
  $
    \calL_1(g^a)^{1/2\ell} 
    \lesssim 
    \calJ(g^a) 
    \lesssim 
    \calL_1(g^a)^{1/(1+2\ell)}.
  $
\end{proposition}
\begin{proof}
  See \Cref{subsec:proof prop:sharpness of rate with poly degree}
\end{proof}

\subsection{Minimizing the surrogate}
\label{subsec:minimizing the surrogate}

Let us now discuss the practical minimization of the surrogate $\calL_1(g)$ defined in \eqref{equ:def L one feature}.
We can first observe that $\calL_1$ is a quadratic function as we can write $\calL_1(g) = \calB_u(g,g)$ where $\calB_u: (h_1, h_2) \mapsto \Expe{ \|\nabla u(\bfX)\|_2^2 \innerp{ \nabla h_1(\bfX)}{\Piperp_{\nabla u(\bfX)} \nabla h_2(\bfX)}_2}$ is bilinear symmetric and nonnegative.
In particular $\calL_1$ is convex, which is to compare with the functional $\calJ$ from \eqref{equ:def of J} based on Poincaré inequalities, which is generally non-convex.

Now assume that $\calG_1$ lies in a finite dimensional vector space $\spanv{\Phi_1, \cdots, \Phi_K}$,  as  in \eqref{equ:Gm vector space} for $m=1$. An element $g\in\calG_1$ admits the representation 
$g(\bfx) = G^T \Phi(\bfx)$ with some   $G\in\bbR^K$. Then $\calL_1(G^T\Phi)$ has a  simple expression in terms of $G\in\bbR^K$.
Indeed, it can be expressed as a positive semi-definite quadratic form of $G$ with some positive semi-definite matrix $H$ which depends on $u$ and $\Phi$, as stated in \Cref{prop:L1 is quadratic form} below.

\begin{proposition}
\label{prop:L1 is quadratic form}
  For any $G\in\bbR^K$ it holds
  \begin{equation}
  \label{equ:L1 is quadratic form}
    \calL_1(G^T \Phi) = G^T H G,
  \end{equation}
  where $H := H^{(1)} - H^{(2)} \in \bbR^{K\times K}$ is a positive semi-definite matrix, with
  \begin{equation}
  \label{equ:def H matrix for L}
  \begin{aligned}
    H^{(1)} &:= \Expe{
      \|\nabla u(\bfX)\|_2^2
      \nabla \Phi(\bfX)^T \nabla \Phi(\bfX)
      }, \\
    H^{(2)} &:= \Expe{
      \nabla \Phi(\bfX)^T \nabla u(\bfX) \nabla u(\bfX)^T \nabla \Phi(\bfX)} .
  \end{aligned}
  \end{equation}
\end{proposition}
\begin{proof}
  Let $G\in\bbR^k$ and $g=G^T\Phi$.
  Using $\|\Piperp_{\nabla u(\bfX)} \nabla g(\bfX) \|_2^2 =  \nabla g(\bfX)^T \Piperp_{\nabla u(\bfX)} \nabla g(\bfX) $ and $\nabla g(\bfX) = \nabla \Phi(\bfX) G$, as well as the linearity of the expectation, we obtain
  \[
    \calL_1(g)
    = G^T \Expe{
      \|\nabla u(\bfX)\|^2
      \nabla \Phi(\bfX)^T \Piperp_{\nabla u(\bfX)} \nabla \Phi(\bfX)
    } G.
  \]
  Then using $\Piperp_{\nabla u(\bfX)} = I_d - \Pi_{\nabla u(\bfX)}$ and $\|\nabla u(\bfX)\|^2 \Pi_{\nabla u(\bfX)} = \nabla u(\bfX) \nabla u(\bfX)^T$ we obtain
  \[
    \calL_1(g)
    = G^T \Expe{
      \nabla \Phi(\bfX)^T 
      \big(
        \|\nabla u(\bfX)\|^2 I_d - \nabla u(\bfX) \nabla u(\bfX)^T
      \big)  
      \nabla \Phi(\bfX)
    } G = G^T H G.
  \]
  The matrix $H\in\bbR^{K\times K}$ is symmetric as the sum of two symmetric matrices and, since $\calL_1(g)\geq 0$ by definition, it is positive semi-definite.
\end{proof}

Now in order to possibly ensure that $\underline{\nu}_{\calK_1,p_1} >0$ and $\overline{\nu}_{\calK_1,p_1} < +\infty$, that is required to properly control the suboptimality of $\calL_1$, we need to restrict the minimization of $\calL_1$ to a compact subset of $\spanv{\Phi_i}_{1\leq i\leq K}$.
A convenient such subset is the unit sphere described in \eqref{equ:Gm vector space sphere}, which we recall in the case $m=1$ is 
\[
  \calG_1 = 
  \Big\{ 
    g: \bfx \mapsto G^T \Phi(\bfx) : 
      G \in \bbR^{K},
      G^T R G = 1
  \Big\},
\]
with $R \in \bbR^{K\times K}$ a symmetric positive definite matrix.
Recall that this restriction to some unit sphere does not modify the minimization problem on $\calJ$.
Indeed since $\calJ(G^T \Phi)$ only depends on $\spanv{G}$, it is equivalent to minimize $\calJ$ over $\spanv{\Phi_1, \cdots, \Phi_K}$ or $\calG_1$, regardless of the choice of the symmetric positive-definite matrix $R$ discussed further below.
In addition to that, minimizing $\calL_1$ over $\calG_1$ is equivalent to finding the minimal generalized eigenpair of the pair $(H, R)$, as we state in \Cref{prop:L1 as generalized ev}.

\begin{proposition}
\label{prop:L1 as generalized ev}
  The minimizers of $\calL_1$ over $\calG_1$ defined in \eqref{equ:Gm vector space sphere} are the functions of the form $g^*(\bfx) = G^{*T} \Phi(\bfx)$, where $G^*$ is a solution to the generalized eigenvalue problem
  \begin{equation}
  \label{equ:L1 as generalized ev}
    \min_{\substack{G \in \bbR^{K} \\G^T R G = 1}} G^T H G,
  \end{equation}
  with $H$  defined in \eqref{equ:def H matrix for L}.
\end{proposition}

Let us briefly discuss on the choice of the matrix $R$.
In \cite{bigoniNonlinearDimensionReduction2022} the authors imposed the condition $\Expe{g(\bfX)^2} = 1$ in order to stabilize the algorithm for minimizing $\calJ$.
This is actually equivalent to taking $R=\Expe{\Phi(\bfX) \Phi(\bfX)^T}$.
However, with our approach, the suboptimality of $\calL_1$ relies on the control of $\|\nabla g(\bfX)\|^2_2$ and not $|g(\bfX)|^2$. 
For example if $R = \Expe{\nabla \Phi(\bfX)^T \nabla \Phi(\bfX)}$, which is equivalent to $\Expe{\|\nabla g(\bfX)\|_2^2} = 1$ for all $g\in\calG_1$, then $\underline{\nu}_{\calK_1,1} = \overline{\nu}_{\calK_1,1}=1$ as defined in \eqref{equ:def nu}, with $\calK_1$ defined in \eqref{equ:def set norm grad features}, that simplifies \Cref{prop:bound J by L,prop:bound L by J,prop:sub optimality} with $p_1=1$.
Hence, when minimizing $\calL_1$, $R = \Expe{\nabla \Phi(\bfX)^T \nabla \Phi(\bfX)}$ appears as a natural choice.

\begin{remark}
  In \cite{romorKernelbasedActiveSubspaces2022}, the authors proposed to select the coefficient matrix $G\in\bbR^{K\times m}$ whose columns are the eigenvector associated to the largest eigenvalues of
  \[
    \tilde{H} = \Expe{\nabla \Phi(\bfX)^+ \nabla u(\bfX) \nabla u(\bfX)^T\nabla \Phi(\bfX)^{+\,T}} \in \bbR^{K\times K},
  \]
  where $A^+$ denotes the Moore-Penrose pseudo inverse of a matrix $A$ of any size.
  It turns out that it can be seen as a special case of our approach.
  First, note that in \Cref{def:poincare inequality}, \Cref{prop:poincare based bound} and \eqref{equ:def L one feature} we decided to use the ambient Euclidean metric of $\bbR^d$ as Riemannian metric on $\calX \subset \bbR^d$.
  However, we could also formulate it equivalently with any other metric $M : \calX \rightarrow \bbR^{d\times d}$, leading to corresponding loss function $\calJ^{M}$, convex surrogate $\calL_1^{M}$ and Poincar\'e constant $C^{M}(\bfX | \calG_1)$.
  Note that this is somehow what is done in \cite{verdiereDiffeomorphismbasedFeatureLearning2025} where the authors leverage what they call Poincar\'e inequalities in feature space.
  Those quantities now involve $\|\cdot\|_{M(\bfX)}^2 = \innerp{\cdot}{M(\bfX) \cdot}_2$ instead of $\|\cdot \|_2^2$ as well as the orthogonal projection matrices  $\Pi^{M(\bfX)}_{\nabla g(\bfX)}$ and $\Pi^{M(\bfX)}_{\nabla u(\bfX)}$ with respect to the inner product matrix $M(\bfX)$.
  In particular, choosing $M(\bfx) = (\nabla \Phi(\bfx) \nabla \Phi(\bfx)^T)^{-1}$ yields
  \[
    \calL_1^{M}(G^T \Phi) = 
    G^T \big(
    \Expe{\|\nabla \Phi(\bfX)^+ \nabla u(\bfX)\|_2^2 \nabla \Phi(\bfX)^T \nabla \Phi(\bfX)^{+\,T}}
    - \tilde{H}
    \big) G,
  \]
  Now if we take $R=\Expe{\|\nabla \Phi(\bfX)^+ \nabla u(\bfX)\|_2^2 \nabla \Phi(\bfX)^T \nabla \Phi(\bfX)^{+\,T}}$ then we have that
  \[
    \min_{\substack{G \in \bbR^{K} \\G^T R G = 1}}
    \calL^M_1(G^T \Phi)
    = \min_{\substack{G \in \bbR^{K} \\G^T R G = 1}}
    (1 - G^T \tilde{H} G)
    = 1 - \max_{\substack{G \in \bbR^{K} \\G^T R G = 1}}
    G^T \tilde{H} G.
  \]
  As a result, the minimizer of $\calL^M_1$ over $\calG_1$ is obtained by finding the eigenvector associated to the largest generalized eigenvalue of $(\tilde{H}, R)$.
  Finally, if by any chance $\Phi$ and $u$ are such that $R=I_K$, then this eigenvector is also solution to \eqref{equ:L1 as generalized ev}.
  In this very specific case, our approach is then equivalent to the one from \cite{romorKernelbasedActiveSubspaces2022}.
\end{remark}

Finally, let us end this section by discussing on the minimization of $\calL_1$ on more general sets $\calG_1$ of features maps.
A general setting is when $\calG_1 = \{g^{(\theta)}: \theta \in \Theta\}$ for some parametrization map $\theta \mapsto g^{(\theta)}$, as it is the case when considering invertible neural networks for feature maps based on diffeomorphism as in \eqref{equ:Gm diffeo}.
With such feature maps, the minimization of $\calL_1$ can only rely on some iterative procedure.
Hence, there might be less incentive to first minimizing $\calL_1$ instead of directly minimizing $\calJ$.
Still, the fact that $\calL_1$ is convex allows to leverage (local) convergence results from convex optimization, and from \Cref{prop:sub optimality bi lip} we may expect $\calL_1$ to be a good approximation of $\calJ$.
Another setting could be to consider feature maps based on low-rank tensor formats and to minimize $\calL_1$ using classical optimization methods such as alternating least squares, for which local convergence results can be found in \cite[Section 9.5.2]{hackbuschTensorSpacesNumerical2019} and which requires at each iteration to solve a quadratic convex optimization problem.
These settings are left to further investigation.

\section{The case of multiple features}
\label{sec:multiple features}

In this section, we propose an extension of the surrogate we introduced in \eqref{equ:def L one feature} to the general case $m> 1$.
The main problem when trying to do so is that for $m>1$, $\spanv{\nabla u(\bfX)}$ and $\spanv{\nabla g(\bfX)}$ are no longer of the same dimension, thus we can no longer expect to relate $\|\Piperp_{\nabla g(\bfX)} \nabla u(\bfx)\|_2$ with $\|\Piperp_{\nabla u(\bfx)} \nabla g(\bfX)\|_F$ similarly as in \Cref{lem:norm projection}.
For this reason, the extension we propose relies on learning one feature at a time, in a coordinate descent or greedy approach.

To do so, let us first introduce a splitting of coordinates of the feature map $g\in\calG_m$.
For $1 \leq j \leq m$ and $g\in\calG_m$ we denote by $g_j$ and $g_{-j}$ respectively the $j$-th and the other $m-1$ components of $g$, and we define the subset of functions which contains exactly the features maps of $\calG_m$ whose $i$-th component is $g_i$ for all $i \neq j$,
\begin{equation}
\label{equ:def Gjmg}
  \calG_{m,g}^{(j)}
  := \{ h \in \calG_m: h_{-j} = g_{-j} \}
  \subset \calG_m.
\end{equation}
The first thing to notice is that for any matrix $W\in\bbR^{d\times (m-1)}$ and vectors $\bfw,\bfv\in\bbR^d$, we can write 
\begin{equation}
\label{equ:split Piperp}
  \Piperp_{(W,\bfw)} \bfv
  = \Piperp_{(W , \Piperp_W \bfw)} \bfv
  = \Piperp_{\Piperp_W \bfw}\Piperp_W \bfv,
\end{equation}
where $(W,\bfw)\in\bbR^{d \times m}$ denotes here the concatenation of $W$ and $\bfw$.
Hence, by considering for any $g\in\calG_m$
\begin{equation}
\label{equ:def w, v}
  w_{g,j}(\bfx) := \Piperp_{\nabla g_{-j}(\bfx)}\nabla g_j(\bfx) \in \bbR^d
  \quad \text{and} \quad
  v_{g,j}(\bfx) := \Piperp_{\nabla g_{-j}(\bfx)} \nabla u(\bfx)\in \bbR^d,
\end{equation}
where $w_{g,j}$ is linear in $g_j$ and $v_{g,j}$ only depends on $g_{-j}$, we can write
\begin{equation}
\label{equ:J with v}
  \calJ(g) = 
   \Expe{\|\Piperp_{w_{g,j}(\bfX)} v_{g,j}(\bfX) \|_2^2
   }.
\end{equation}
We now end up with an expression for $\calJ(g)$ which is very similar to what we had in the case $m=1$, where $\nabla g$ and $\nabla u$ are replaced respectively by $w_{g,j}$ and $v_{g,j}$.
For this reason, we introduce the following surrogate,
\begin{equation}
\label{equ:def of Lm}
  \calL_{m,j}(g) 
  := \Expe{
    \|v_{g,j}(\bfX)\|_2^2 \|\Piperp_{v_{g,j}(\bfX)} w_{g,j}(\bfX) \|_2^2
  },
\end{equation}
with $w_{g,j}, v_{g,j}$ as defined in \eqref{equ:def w, v}.
Now if $\calG_m$ satisfies \Cref{assump:full rank jacobian as}, then $\|w_{g,j}(\bfX)\|_2 > 0$ almost surely, thus taking $w=w_{g,j}(\bfX)$ and $v=v_{g,j}(\bfX)$ in \Cref{lem:norm projection} yields an expression similar to \eqref{equ:norm projection applied 1},
\begin{equation}
\label{equ:norm projection m}
\begin{aligned}
  \calJ(g)
  &= \Expe{
    \|w_{g,j}(\bfX)\|_2^{-2}
    \|v_{g,j}(\bfX)\|_2^2
    \|\Piperp_{v_{g,j}(\bfX)} w_{g,j}(\bfX)\|_2^2
  }, \\
  \calL_{m,j}(g)
  &= \Expe{
    \|w_{g,j}(\bfX)\|_2^2
    \|\Piperp_{w_{g,j}(\bfX)} v_{g,j}(\bfX)\|_2^2
  }.
\end{aligned}
\end{equation}
In the first and second right-hand sides of \eqref{equ:norm projection m} we recognize the terms within the expectation in respectively \eqref{equ:J with v} and \eqref{equ:def of Lm}, up to a factor $\|w_{g,j}(\bfX)\|_2^2$.
Note that $\calG_m$ satisfies \Cref{assump:full rank jacobian as} if and only if for all $g\in\calG_m$ and $1\leq j\leq m$, $\|w_{g,j}(\bfX)\|_2 > 0$ almost surely.

The structure of this section is as follows.
First in \Cref{subsec:sub optimality m} we discuss on the suboptimality of $\calL_{m,j}$.
Then in \Cref{subsec:minimizing the surrogate m} we discuss on the problem of minimizing $\calL_{m,j}$.

\subsection{suboptimality}
\label{subsec:sub optimality m}

Similarly to \Cref{prop:exact recovery}, we can show that in the case $u=f\circ g$ for some $g\in\calG_m$, minimizing $\calL_{m,j}$ is equivalent to minimizing $\calJ$, as stated in \Cref{prop:exact recovery m} below.

\begin{proposition}
\label{prop:exact recovery m}
  Assume that $\calG_m$ satisfies \Cref{assump:full rank jacobian as}, then for all $1\leq j\leq m$ and all $g\in\calG_m$ we have $\calJ(g)=0$ if and only if $\calL_{m,j}(g) = 0$.
\end{proposition}
\begin{proof}
  Let $g\in\calG_m$ and $1\leq j\leq m$.
  Firstly, assume that $\calJ (g) = 0$.
  Then using \eqref{equ:J with v} yields that almost surely $\|\Piperp_{w_{g,j}(\bfX)} v_{g,j}(\bfX) \|_2^2 = 0$, thus \eqref{equ:norm projection m} yields $\calL_{m,j}(g)=0$.
  Secondly, assume that $\calL_{m,j}(g) = 0$.
  Then \eqref{equ:norm projection m} implies that almost surely, $\|v_{g,j}(\bfX)\|_2^2=0$ or $\|\Piperp_{v_{g,j}(\bfX)} w_{g,j}(\bfX) \|_2^2=0$.
  If $\|v_{g,j}(\bfX)\|_2^2=0$ then $\nabla u(\bfX) \in \spanv{\nabla g_{-j}(\bfX)} \subset \spanv{\nabla g(\bfX)}$, thus $\calJ(g) = 0$.
  Otherwise, $\|\Piperp_{v_{g,j}(\bfX)} w_{g,j}(\bfX) \|_2^2 = 0$, thus using \Cref{lem:norm projection}, using that almost surely $\|w_{g,j}(\bfX)\|_2 >0$ from \Cref{assump:full rank jacobian as}, and using \eqref{equ:split Piperp} with $W=\nabla g_{-j}(\bfX)$, $\bfw = w_{g,j(\bfX)}$ and $\bfv = v_{g,j}(\bfX)$, we obtain almost surely
  \[
    \Piperp_{w_{g,j}(\bfX)} v_{g,j}(\bfX) 
    = \Piperp_{\Piperp_{\nabla g_{-j}(\bfX)} \nabla g_j(\bfX)} 
    \Piperp_{\nabla g_{-j}(\bfX)} \nabla u(\bfX)
    =\Piperp_{\nabla g(\bfX)} \nabla u(\bfX)
    = 0,
  \]
  which implies $\calJ(g)=0$ and concludes the proof.
\end{proof}

Now in general we cannot expect the minimum of $\calJ$ to be zero, and we would want to quantify the suboptimality of $\calL_{m,j}$.
In view of \eqref{equ:norm projection m}, we then need some control on $\|w_{g,j}(\bfX)\|_2^2$ defined in \eqref{equ:def w, v}.
The new difficulty compared to the case $m=1$ from \Cref{sec:one feature} lies in lower bounding this quantity.
A first relevant result is given in the following \Cref{lem:bounds on norm w}.

\begin{lemma}
\label{lem:bounds on norm w}
  For all $1\leq j\leq m$, all $g \in \calC^1(\calX, \bbR^m)$ and all $\bfx \in \calX$ it holds
  \begin{equation*}
    \sigma_{m}(\nabla g(\bfx))^2
    \leq \|w_{g,j}(\bfx)\|_2^2
    \leq \|\nabla g_j(\bfx)\|_2^2,
  \end{equation*}
  with $w_{g,j}$ defined in \eqref{equ:def w, v}.
\end{lemma}
\begin{proof}
  See \Cref{subsec:proof lem:bounds on norm w}.
\end{proof}

In view of the above \Cref{lem:bounds on norm w},
the ideal setting to quantify suboptimality of $\calL_{m,j}$ would be the bi-Lipschitz setting as in \Cref{prop:sub optimality bi lip}, which we restate in the context $m\geq 1$ in \Cref{prop:sub optimality bi lip Lm}. 
Note again that such context occurs when considering $\calG_m$ as in \eqref{equ:Gm diffeo} with $\calD$ containing only bi-Lipschitz diffeomorphisms whose associated constants are uniformly bounded over $\calD$.
Note also that our analysis in \Cref{prop:sub optimality bi lip Lm} is rather coarse, and with finer assumptions on bi-Lipschitz constant depending on the $j$-th feature, one could probably obtain a sharper bound.

\begin{proposition}
\label{prop:sub optimality bi lip Lm}
  Assume that there exist $c>0$ and $C < +\infty$ such that for all $g\in\calG_m$ and $\bfx\in\calX$, $c \leq \sigma_{m}(\nabla g(\bfx))^2 \leq \sigma_{1}(\nabla g(\bfx))^2 \leq C$.
  Then for all $1\leq j\leq m$ and $g\in\calG_m$,
  \begin{equation}
    c \calJ(g)
    \leq 
    \calL_{m,j}(g) 
    \leq
    C \calJ(g).
  \end{equation}
  Moreover for all $1\leq j\leq m$ and all $g \in \calG_m$, if $g^*$ is a minimizer of $\calL_{m,j}$ over $\calG_{m,g}^{(j)}$ defined in \eqref{equ:def Gjmg}, then
  \begin{equation}
    \calJ(g^*) 
    \leq (C / c) 
    \inf_{h\in \calG_{m,g}^{(j)}} \calJ(h)
  \end{equation}
\end{proposition}
\begin{proof}
  Let $g\in\calG_m$ and $1\leq j\leq m$.
  Using \Cref{lem:bounds on norm w} and $\|\nabla g_j(\bfX)\|_2^2 \leq \sigma_1(\nabla g(\bfX))^2$, we obtain $c \leq \|w_{g,j}(\bfX)\|_2^2 \leq C$.
  Then using \eqref{equ:J with v} and \eqref{equ:norm projection m} we obtain $c\calJ(g) \leq \calL_{m,j}(g) \leq C \calJ(g)$.
  The suboptimality result then directly follows.
\end{proof}

Without the bi-Lipschitz assumption, the suboptimality of our surrogate should rely on controlling the probability of small and large deviations of $\|w_{g,j}(\bfX)\|_2^2$, as in the setting $m=1$.
For the large deviations, it is rather straightforward since we have $\|w_{g,j}(\bfX)\|_2^2 \leq \|\nabla g_j(\bfX)\|_2^2$, thus we can re-use the results from \Cref{sec:Large and small deviations inequalities}.
However for the small deviations, meaning controlling the probability of the event $(\|w_{g,j}(\bfX)\|^2_2 \leq \alpha)$ for $\alpha >0$, the task is much more complicated.
The main issue is that the family $\{\nabla g_i(\bfX)\}_{1\leq i\leq m}$ is not orthogonal in general.

Still, we can obtain small deviation inequalities, and thus establish a suboptimality result, for the specific case of multiple polynomial features with  bounded total degree, as stated in \Cref{prop:suboptimality uniform m} below.
This can be done by combining \Cref{prop:small deviation sigma m} and \Cref{lem:bounds on norm w}, which will involve
\begin{equation}
\label{equ:def Km}
  \calK_m := \{\bfx \mapsto \det(\nabla g(\bfx)^T \nabla g(\bfx)) : g\in\calG_m\}.
\end{equation}
The main idea behind our result is that if $\calG_m$ only contains polynomials of total degree at most $\ell+1$, then $\calK_m$ only contains polynomials of total degree at most $2\ell m$, thus elements of $\calK_m$ satisfy \eqref{equ:remez inequality multivariate}.
Note that it remains unclear whether other classes of feature maps satisfy such property, hence we restrict our analysis to polynomials.
Note also that we state our result for features such that $\Expe{\|\nabla g(\bfX)\|_2^2} \leq 1$ for all $1\leq j\leq m$.
Similarly to the assumption on $u$, this is only for readability purpose, as results without this assumption can be obtained by normalizing the feature maps and leveraging the properties of $\calL_{m,j}$ and $\calJ$.
Moreover, such assumption on $\calG_m$ can be ensured by taking $\calG_m$ as in \eqref{equ:Gm vector space sphere} with $\Phi :\bbR^d \rightarrow \bbR^K$ a polynomial and $R = \Expe{\nabla \Phi(\bfX)^T \nabla \Phi(\bfX)}$.

\begin{proposition}
\label{prop:suboptimality uniform m}
  Assume that $\bfX$ is an absolutely continuous random variable on $\bbR^d$ whose distribution is $s$-concave with $s\in(0,1/d]$.
  Assume that every $g\in\calG_m$ is a non-constant polynomial with total degree at most $\ell+1 \geq 2$ such that $\Expe{\|\nabla g_j(\bfX)\|_2^2} \leq 1$ for all $1\leq j\leq m$.
  Assume that $\sup_{\bfx\in\calX} \|\nabla u(\bfx)\|_2^2 \leq 1$.
  Then for all $1\leq j \leq m$, $g\in\calG_m$ and $p_1\geq 1$ it holds
  \[
    \tilde{\gamma_2}^{-1} \calL_{m,j}(g)
    \leq 
    \calJ(g)
    \leq 
    \tilde \gamma_1
    \underline{\nu}_{\calK_m, p_1}^{-\frac{1}{1+2\ell m}}
    \calL_{m,j}(g)^{\frac{1}{1 + 2\ell m}},
  \]
  with $\underline{\nu}_{\calK_m, p_1}$ defined in \eqref{equ:def nu}, with $\calK_m$ defined in \eqref{equ:def Km}, and with $\tilde \gamma_1$ and $\tilde \gamma_2$ both defined in \Cref{subsec:proof prop:suboptimality uniform m} satisfying
  \[
    \tilde \gamma_1 \leq 2^9 m ^{\frac{1}{4\ell}} s^{-1} \min\{s^{-1},3\ell p_1 m\},
    \quad 
    \tilde \gamma_2 \leq 2^{1+6\ell} s^{-2\ell}.
  \]
  Moreover for all $g \in \calG_m$, if $g^*$ is a minimizer of $\calL_{m,j}$ over $\calG_{m,g}^{(j)}$, it holds
  \begin{equation}
    \calJ(g^*)  
    \leq 
    \tilde{\gamma}_3 
    \underline{\nu}_{\calK_m, p_1}^{-\frac{1}{1+2\ell m}}
    \inf_{h\in \calG_{m,g}^{(j)}} \calJ(h)^{\frac{1}{1+2\ell m}},
  \end{equation}
  with $\tilde \gamma_3$ defined in \Cref{subsec:proof prop:suboptimality uniform m} and satisfying $\tilde \gamma_3 \leq 2^{10} m^{\frac{1}{4\ell}} s^{-1} \min\{s^{-1}, 3\ell p_1 m\}$.
\end{proposition}
\begin{proof}
  See \Cref{subsec:proof prop:suboptimality uniform m}.
\end{proof}

Let us end this section by pointing out the main differences between the suboptimality result for $m=1$ feature in \Cref{cor:suboptimality uniform} and the above suboptimality result for $m>1$ features in \Cref{prop:suboptimality uniform m}.
The first difference is that for multiple features, the result actually concerns only one of the $m$ features.
The second difference is that the exponent is now $1/(1+2\ell m)$ instead of $1/(1+ 2\ell)$, which is much worse. 
In particular, the larger is $m$, in other words the more features we have to learn, the worse is the suboptimality result.
The third difference is that we do not know yet how to properly control $\underline{\nu}_{\calK_m, p_1}$.

\subsection{Minimizing the surrogate}
\label{subsec:minimizing the surrogate m}

Now we discuss on the practical minimization of the surrogate  $\calL_{m,j}$ introduced in \eqref{equ:def of Lm}.
At first glance, $\calL_{m,j}$ does not have much better properties than $\calJ$ since it involves $\Pi_{\nabla g_{-j}(\bfX)}$.
However, optimizing the $j$-th coordinate is actually similar to the one feature case.
Indeed for all $1\leq j\leq m$ and $g \in \calG_m$, $\calL_{m,j}$ is quadratic on $\calG_{m,g}^{(j)}$, more precisely
\[
  h \mapsto \calL_{m,j}((g_1, \cdots, g_{j-1}, h, g_{j+1}, \cdots, g_m))
\]
is quadratic. 
Indeed, we can write $\calL_{m,j}(g) = \calB_{u,g_{-j}}(g_j,g_j)$ where $\calB_{u,g_{-j}}: (h_1, h_2) \mapsto\Expe{\|v_{g,j}(\bfX)\|_2^2 \innerp{ \nabla h_1(\bfX)}{\Piperp_{v_{g,j}(\bfX)} \Piperp_{\nabla g_{-j}(\bfX)} \nabla h_2(\bfX)}_2}$ is bilinear symmetric and nonnegative.
Indeed, $\Piperp_{v_{g,j}(\bfX)} $ and $\Piperp_{\nabla g_{-j}(\bfX)}$ commute since $v_{g,j}(\bfX) \in \spanv{\nabla g_{-j}(\bfX)}^{\perp}$.
In particular $\calL_{m,j}$ is convex on $\calG_{m,g}^{(j)}$, which is to compare with $\calJ$ which is in general non-convex on $\calG_{m,g}^{(j)}$.
Now again when $\calG_{m,g}^{(j)} \subset \spanv{\Phi_i}_{1\leq i\leq K}$ we can express 
\[
  G_j \mapsto \calL_{m,j}((g_1, \cdots, g_{j-1}, G_j^T \Phi, g_{j+1}, \cdots, g_m))
\]
as a positive semi-definite quadratic form on $G_j$ with some positive semi-definite matrix $H_{g,j}$ which depends on $u$, $\Phi$ and $g_{-j}$, as stated in \Cref{prop:Lm is quadratic form} below.

\begin{proposition}
\label{prop:Lm is quadratic form}
  For all $g\in\calG_m$ and $G_j\in\bbR^K$ such that $h\in\calG_{m,g}^{(j)}$ and $h_j = G_j^T \Phi$, it holds
  \begin{equation}
  \label{equ:Lm is quadratic form}
    \calL_{m,j}(h) = G_j^T H_{g,j} G_j,
  \end{equation}
  where $H_{g,j} = H_{g,j}^{(1)} - H_{g,j}^{(2)} \in \bbR^{K\times K}$ is a positive semi-definite matrix with
  \begin{equation}
  \label{equ:def Hm matrix for Lm}
  \begin{aligned}
    H_{g,j}^{(1)} &:= \Expe{
      \|\Piperp_{\nabla g_{-j}(\bfX)} \nabla u(\bfX)\|_2^2
      \nabla \Phi(\bfX)^T 
      \Piperp_{\nabla g_{-j}(\bfX)}
      \nabla \Phi(\bfX)
      }, \\
    H_{g,j}^{(2)} &:= \Expe{
      \|\Piperp_{\nabla g_{-j}(\bfX)} \nabla u(\bfX)\|_2^2
      \nabla \Phi(\bfX)^T 
      \nabla u(\bfX) \Piperp_{\nabla g_{-j}(\bfX)} \nabla u(\bfX)^T
      \nabla \Phi(\bfX)
      },
    \end{aligned}
  \end{equation}
  which depend on $u$, $\Phi$ and $g_{-j}$.
\end{proposition}
\begin{proof}
  See \Cref{subsec:proof prop:Lm is quadratic form}.
\end{proof}

Now, let us study the case when $\calG_m$ has the form \eqref{equ:Gm vector space sphere}.
In this case, we have a constrained minimization of a Rayleigh quotient, where the constraints are $G_i^T R G_j = \delta_{i,j}$.
It turns out that this minimization problem is equivalent to some generalized eigenvalue problem, with orthogonality matrix $R$.
Indeed, by noticing that $H_{g,j}G_{-j} = 0$, we have for any $\alpha$ that $H_{g,j} + \alpha R G_{-j} G_{-j}^T R$ admits $G_{-j}$ as generalized eigenvector with associated eigenvalue $\alpha$.
Now to be sure that this is not the smallest one, we shall take $\alpha$ large enough, in particular larger than any generalized eigenvalue $\lambda$ of $(H_{g,j}, R)$.
Note that by definition of $\calL_{m,j}$ in \eqref{equ:def of Lm}, it holds
\[
  G_j H_{g,j} G_j 
  \leq \Expe{\|\nabla u(\bfX)\|_2^2 \|\nabla g_j(\bfX)\|_2^2}
  = G_j^T H^{(1)} G_j,
\]
with $H^{(1)}$ defined in \eqref{equ:def H matrix for L}.
Hence, choosing $\alpha$ as the largest generalized eigenvalue of $(H^{(1)}, R)$ yields that $G_{-j}$ is not associated to the smallest generalized eigenvector of $(H_{g,j} + \alpha R G_{-j} G_{-j}^T R, R)$, and that the corresponding generalized eigenvector is orthogonal to $G_{-j}$ with respect to $R$.
This is summarized in \Cref{prop:Lm as generalized ev} below.

\begin{proposition}
\label{prop:Lm as generalized ev}
  Let $g = G^T \Phi \in \calG_m$ as defined in \eqref{equ:Gm vector space sphere}.
  The minimizers of $h \mapsto \calL_{m,j}(h)$ over $\calG_{m,g}^{(j)}$ defined in \eqref{equ:def Gjmg} are the functions of the form $g^* = G^{*\,T}\Phi$, where $G_i^* = G_i$ for all $i\neq j$ and $G_j^*$ is a solution to the generalized eigenvalue problem
    \begin{equation}
    \label{equ:Lm as generalized ev}
      \min_{\substack{v \in \bbR^{K} \\v^T R v = 1}}
      v^T (H_{g,j} + \alpha H_{g,j}^{(3)}) v
    \end{equation}
    with $H_{g,j}$ as defined in \eqref{equ:def Hm matrix for Lm}, $H_{g,j}^{(3)} := R G_{-j} G_{-j}^T R$ and $\alpha$ as prescribed above.
\end{proposition}

With \Cref{prop:Lm as generalized ev}, we have an explicit formulation to compute the $j$-th feature $g_j = G_j^T \Phi$ from the $j-1$ previously learnt features $g_{-j} = G_{-j}^T \Phi$, based on the surrogate we introduced in \eqref{equ:def of Lm}.
We now propose to leverage it in a greedy algorithm, in which   we compute successively $G_1, \cdots, G_m$.
The algorithm is described in \Cref{algo:greedy}.
Note that in practice, at the end of each greedy iteration, one may consider running a few iterations of some minimization algorithm on $G \mapsto \calJ(G^T \Phi)$.
Indeed, one shall keep in mind that the purpose of our surrogates $\calL_1$ and $\calL_{m,j}$ is to obtain a suboptimal solution to the minimization of $\calJ$.
Although the latter optimization problem is rather difficult to solve, we can still leverage first and second order optimization algorithms, as proposed in \cite{bigoniNonlinearDimensionReduction2022}.

\begin{algorithm}[H]
\caption{Greedy algorithm to learn multiple features}
\label{algo:greedy}
\begin{algorithmic}[1]
  \Require Number of features $m\geq 2$, basis $\Phi : \bbR^d \rightarrow \bbR^K$ of feature space, inner product matrix $R\in\bbR^{K\times K}$.
  \Ensure $G\in\bbR^{k \times m}$ and $G^T R G = I_m$.
  \State Compute $H^{(1)},H^{(2)} \in \bbR^{K\times K}$ defined in \eqref{equ:def H matrix for L}.
  \State Compute $G_1 \in \bbR^K$ solution to \eqref{equ:L1 as generalized ev} and $\alpha$ largest generalized eigenvalue of $(H^{(1)}, R)$.
  \State Set $G \gets (G_1, 0, \cdots, 0) \in \bbR^{K\times m}$ and $g=G^T \Phi$.
  \For{$2\leq j\leq m$}
    \State Compute $H_{g,j}^{(1)}, H_{g,j}^{(2)},H_{g,j}^{(3)} \in \bbR^{K\times K}$ from \eqref{equ:def Hm matrix for Lm} and \eqref{equ:Lm as generalized ev}.
    \State Compute $G_j^*$ solution to \eqref{equ:Lm as generalized ev}.
    \State Set $G \gets (G_1, \cdots , G_{j-1} , G_j^*, 0, \cdots, 0) \in \bbR^{K\times m}$.
  \EndFor
  \end{algorithmic}
\end{algorithm}

In the above \Cref{algo:greedy} we chose to use our surrogates $(\calL_{m,j})_{1\leq j\leq m}$ for learning (or initializing) the coefficients $G$ of a feature map in a single pass.
An alternative approach could be to do multiple passes with our surrogates, in a coordinate-descent manner.
However, in practice we did not observe significant improvements in the numerical experiments compared to our greedy approach.
For this reason we did not include such procedure in the present work and preferred to keep the simpler greedy approach described in \Cref{algo:greedy}.
Still, there is room for improvement.

\section{Numerical experiments}
\label{sec:numerical experiments}

\subsection{Setting}
\label{subsec:setting}

\paragraph{Benchmark functions.}
We consider the following four benchmark functions $u_i : \bbR^d \rightarrow \bbR$ with $d=8$:
\[
\begin{aligned}
    u_1(\bfx) := \sin(\frac{4}{\pi^2} \bfx^T \bfx),& \quad
    u_2(\bfx) := \cos(\frac{1}{2} \bfx^T \bfx) + \sin(\frac{1}{2} \bfx^T M \bfx), \\
    u_3(\bfx) := \exp(\frac{1}{d} \sum_{i=1}^d \sin(x_i) e^{\cos(x_i)}),& \quad
    u_4(\bfx) := \frac{2 \pi x_3 (x_4 - x_6)}{
        \ln(x_2 / x_1)
        (1 + \frac{2 x_7 x_3}{\ln(x_2 / x_1) x_1^2 x_8}
        + \frac{x_3}{x_5}).
    }
\end{aligned}
\]
In $u_2$, we choose the matrix $M = (\frac{1}{i+j-1})_{1\leq i,j\leq d} \in\bbR^{d\times d}$.
For $u_1,u_2$ and $u_3$ we take $\bfX \sim \calU(]-\frac{\pi}{2},\frac{\pi}{2}[^d)$.
For the so-called Borehole function $u_4$ we take $X_1 \sim \calN(0.1, 0.0161812)$, $X_2 \sim \log \calN(7.71, 1.0056)$, $X_3 \sim \calU(63070, 115600)$, $X_4 \sim \calU(990, 1110)$, $X_5 \sim \calU(63.1, 116)$, $X_6 \sim \calU(700, 820)$, $X_7 \sim \calU(1120, 1680)$ and $X_8 \sim \calU(9855, 12045)$.

The benchmark $u_1$ illustrates \Cref{prop:exact recovery}, as it can be expressed as a univariate function of a polynomial, while not being a polynomial.
Thus taking $\calG_1$ as a suitable set of polynomials, we know that $\inf_{\calG_1} \calJ = 0$, and in that case the minimizers of $\calJ$ and $\calL_1$ are the same.
Hence, $u_1$ is perfectly suited for the surrogate we introduced in \Cref{sec:one feature}.
Similarly, the benchmark $u_2$ is such that $\inf_{\calG_2} \calJ = 0$ for a suitable set of polynomials $\calG_2$.
Hence, $u_2$ seems well suited for the surrogate we introduced in \Cref{sec:multiple features}, although we have less theoretical guarantees.
Finally, the benchmarks $u_3$ and $u_4$ were investigated in the context of nonlinear dimension reduction respectively in \cite{verdiereDiffeomorphismbasedFeatureLearning2025} and \cite{bigoniNonlinearDimensionReduction2022}.

\paragraph{Quantities monitored.}
In our experiments, we will monitor $4$ quantities.
The first two are the Poincar\'e inequality based quantity $\calJ(g)$ defined in \eqref{equ:def of J} and the final approximation error $e_g(f)$ defined by 
\[
    e_g(f) := \|u - f \circ g \|_{L^2}.
\]
We estimate these quantities with their Monte-Carlo estimators on test samples $\Xi^{test} \subset \calX$ of sizes $N^{test}=1000$, not used for learning.
We also monitor the Monte-Carlo estimators $\widehat{\calJ}(g)$ and $\widehat{e}_g(f)$ on some training set $\Xi^{train}\subset\calX$ of various sizes $N^{train}$, which will be the quantities directly minimized to compute $g$ and $f$.
More precisely, we draw $20$ realizations of $\Xi^{train}$ and $\Xi^{test}$ and monitor the quantiles of those $4$ quantities over those $20$ realizations.

\paragraph{Sets of feature maps.}

We consider feature maps of the form \eqref{equ:Gm vector space sphere} with $\Phi : \bbR^d \rightarrow \bbR^K$ a multivariate polynomial obtained by sparse tensorizations of univariate orthonormal polynomials.
Those are obtained by selecting multi-indices with $p\in[0, +\infty]$ norm bounded by $k\in [1, +\infty)$, where we see $(p,k)$ as hyperparameter that we will optimize using a cross-validation procedure.

More precisely, for a product measure $d\mu(\bfx) = \bigotimes_{\nu=1}^d d\mu_{\nu}(x_\nu)$ and $1\leq \nu \leq d$, we consider $(\phi_i^{\nu})_{i\in\bbN}$ an orthonormal polynomial basis of $L^2_{\mu_{\nu}}$, meaning $\Expe{\phi_i^{\nu}(X_{\nu}) \phi_j^{\nu}(X_{\nu})} = \delta_{i,j}$, and $\phi_i^{\nu}$ has degree $i$.
Then $(\phi_{\alpha})_{\alpha \in \bbN^d}$ is an orthonormal basis of $L^2_{\mu_{\nu}}$, where for a multi-index $\alpha \in \bbN^d$ we define the multivariate polynomial $\phi_{\alpha}$ as
\[
    \phi_{\alpha}(\bfx) := \prod_{\nu=1}^{d}\phi_{\alpha_{\nu}}^{\nu}(x_{\nu}).
\]
The coordinates of $\Phi$ are then given by multi-indices from some prescribed set $\Lambda \subset \bbN^d$.
Here we consider sets of multi-indices of the form $\Lambda_{p,k}:= \{\alpha \in \bbN^{d} : \|\alpha\|_p \leq k\} \setminus \{0\}$, which yields an associated basis $\Phi = \Phi^{\Lambda_{p,k}}$.
Note that we excluded the multi-index $0$ so that \Cref{assump:full rank jacobian as} is satisfied.
Note also that such definition ensures that $id \in \spanv{\Phi_1, \cdots, \Phi_K}$, thus $\calG_m$ contains all linear feature maps, including the one corresponding to the active subspace method.
The hyperparameters $(p,k)$ are selected  using a cross-validation procedure described bellow.

\paragraph{Computing feature maps.}

For fixed hyperparameters $(p,k)$, associated basis $\Phi: \bbR^d \rightarrow \bbR^K$ and a fixed training set $\Xi^{train} \subset \calX$ of size $N^{train}$, we aim to minimize the Monte-Carlo estimation $\bbR^{K\times m}\ni G \mapsto \widehat{\calJ}(G^T \Phi)$ evaluated from the samples in $\Xi^{train}$.
To do so, we consider $3$ learning procedures.

The first procedure, which we consider as the reference, is based on a preconditioned nonlinear conjugate gradient algorithm on the Grassmann manifold $Grass(m, K)$.
We use $\hat{\Sigma}(G) \in \bbR^{K\times m}$ as preconditioning matrix at point $G\in\bbR^{K\times m}$, which is the Monte-Carlo estimation of $\Sigma(G)$ defined in \cite[Proposition 3.2]{bigoniNonlinearDimensionReduction2022}.
We choose as initial point the matrix $G^0\in\bbR^{K\times m}$ which minimizes $\widehat{\calJ}$ on the set of linear features, which corresponds to the active subspace method.
We denote this reference procedure as GLI, standing for Grassmann Linear Initialization.
The second procedure is based on \Cref{algo:greedy}, where the surrogates $\calL_{m,j}$ are estimated via Monte-Carlo based on $\Xi^{train}$.
Within \Cref{algo:greedy} we use the orthogonality matrix $R = \Expe{\nabla \Phi(\bfX)^T \nabla \Phi(\bfX)}$.
We denote this procedure as SUR, standing for SURrogate.
The third procedure is the same as the first one, but where the initial point is given by the second procedure. 
We denote this procedure as GSI, standing for Grassmann Surrogate Initialization.
Note that after procedure, we orthonormalize again $g$ such that $\Expe{\nabla g(\bfX)^T \nabla g(\bfX)} = I_m$ in order to obtain better stability of the subsequent regression problem.
Finally, for each procedure, we select specific hyperparameters $(p,k)$ using a $5$-fold cross-validation procedure on $\Xi^{train}$, among the $10$ following values,
\[
    (p,k) \in 
    \big( \{0.8\} \times \{2,3,4,5\} \big) \cup
    \big( \{0.9\} \times \{2,3,4\} \big) \cup
    \big( \{1\} \times \{1,2,3\} \big).
\]
Note that this set of hyperparameters has been chosen arbitrarily to ensure a compromise between computational cost and high polynomial degree.

\paragraph{Regression function.}

Once we have learnt $g\in\calG_m$, we then perform a classical regression task to learn a regression function $f:\bbR^m \rightarrow \bbR$, with $g(\bfX)\in\bbR^m$ as input variable and $u(\bfX)\in\bbR$ as output variable.
In particular here, we have chosen to use kernel ridge regression with Gaussian kernel $\kappa(\bfy, \bfz) := \exp(-\gamma \|\bfy - \bfz\|^2_{2})$ for any $\bfy,\bfz\in\bbR^m$ and some hyperparameter $\gamma>0$.
Then with $\{\bfz^{(i)}\}_{1\leq i\leq N^{train}} := g(\Xi^{train}) \subset \bbR^m$, we consider
\[
    f: \bfz  \mapsto \sum_{i=1}^{N^{train}} a_i \kappa(\bfz^{(i)}, \bfz),      
\]
with $\bfa :=  ( K + \alpha I_N)^{-1} \bfu \in \bbR^{N^{train}}$ for some regularization parameter $\alpha >0$, where $K := (\kappa(\bfz^{(i)}, \bfz^{(j)}))_{1\leq i,j\leq N^{train}}$ and $\bfu := (u(\bfx^{(i)}))_{1\leq i\leq N^{train}}$.
Here the kernel parameter $\gamma$ and the regularization parameter $\alpha$ are hyperparameters learnt using a $10$-fold cross-validation procedure, such that $\log_{10}(\gamma)$ is selected from $30$ points uniformly spaced in $[-6, -2]$, and $\log_{10}(\alpha)$ is selected from $40$ points uniformly spaced in $[-11, -5]$.
Note that these sets of hyperparameters have been chosen arbitrarily to ensure a compromise between computational cost and flexibility of the regression model.

Note that we followed a procedure similar to \cite[Algorithm 4]{bigoniNonlinearDimensionReduction2022}.
In particular the hyperparameters for $g$ and $f$ are fitted via two separate cross-validation procedures.
One may also consider optimizing both in one global cross-validation procedure, although it would be much more costly.

\begin{remark}
    In the present work we chose to use a kernel ridge regression with Gaussian kernel to approximate the conditional expectation $\bfz \mapsto \Expe{u(\bfX)|g(\bfX)=\bfz}$.
    The main motivation for this choice was that it is more stable than polynomial regression, especially in a low training sample size regime.
    Note that with additional information on the smoothness of $\bfz \mapsto \Expe{u(\bfX)|g(\bfX)=\bfz}$ one may consider using a Mat\`ern kernel instead, with adapted smoothness hyperparameter.
    We leave this to further investigation, although it is an important question to better tackle the regression step.
\end{remark}

\paragraph{Implementation and data statement.}

The cross-validation procedures as well as the Kernel ridge regression rely on the library \emph{sklearn} \cite{scikit-learn}.
The optimization on Grassmann manifolds rely on the library \emph{pymanopt} \cite{JMLR:v17:16-177}.
The orthonormal polynomials feature maps rely on the python library \emph{tensap} \cite{anthonynouyAnthonynouyTensapV152023}.
The code underlying this article is freely available in \cite{nouyAlexandrepascoTensapV162poincarelearningpaper2025} at \url{https://doi.org/10.5281/zenodo.15430309}.

\subsection{The case of a single feature}
\label{subsec:results one features}

In this section we investigate the performances of the three learning procedures GLI, SUR and GLS when learning $m=1$ feature to approximate $u_1$, $u_3$ and $u_4$.

\begin{figure}[H]
    \centering
    \includegraphics[page=1, width=0.99\textwidth]{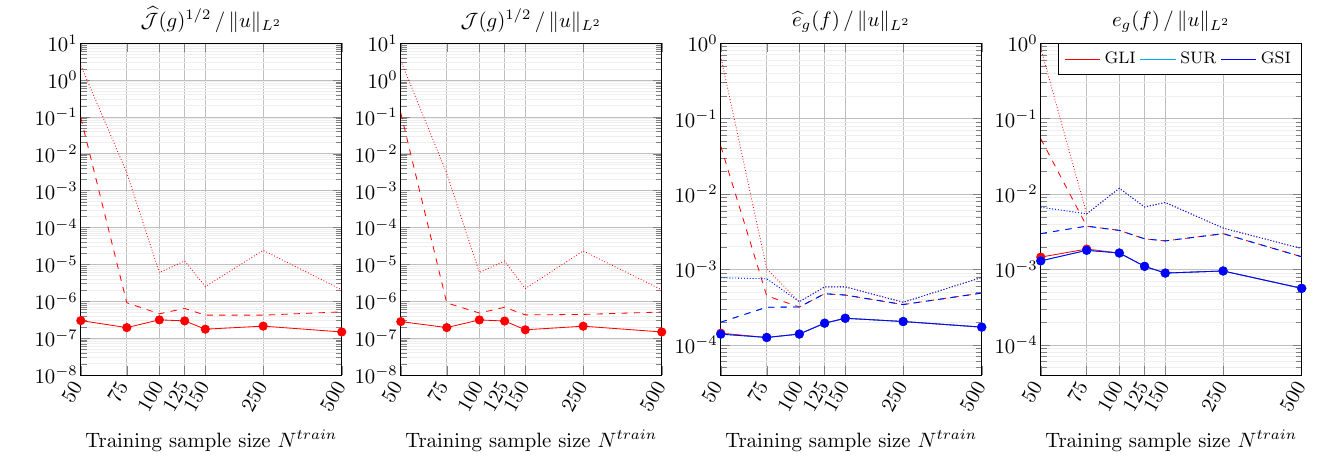}
    \caption[caption]{\footnotesize 
    Evolution of quantiles with respect to the size of the training sample for $u_1$ with $m=1$.
    The quantiles $50\%$, $90\%$ and $100\%$ are represented respectively by the continuous, dashed and dotted lines.
    }
    \label{fig:u1 vs ntrain}
\end{figure}

Let us start with $u_1$ for which results are displayed in \Cref{fig:u1 vs ntrain}.
Recall that this benchmark function can be exactly expressed as a function of $m=1$ polynomial feature of total degree $2$.
Note that, although the graphs in the two left plots seem to coincide exactly, this is not the case, and it is due to the large vertical scale.
This means that in this example the Monte-Carlo estimator $\widehat{\calJ}(g)$ is a good estimator of $\calJ(g)$, even for small values of $N^{train}$.

Firstly for all $N^{train}$, we observe that SUR and GSI always yield the minimizer of $\widehat{\calJ}$, which turns out to be $0$ as for $\calJ$.
This is the reason why they are not displayed in the first two plots, as they are equal to machine zero, and why SUR and GSI yield the same results, as the minimum is already achieved with SUR. 
This could be expected from \Cref{prop:exact recovery}, although it would probably have failed for even smaller $N^{train}$.
As a result, we know that the minimum of $e_g$ over all possible regression functions is $0$.
However, we still observe that the approximation errors $e_g(f)$ may remain rather high, in particular for $u_2$.
This is most probably because kernel ridge regression is not suited for approximating $\bfz \mapsto \Expe{u(\bfX) | g(\bfX) = \bfz}$, especially for small training sets.
In particular, since $\widehat{e}_g(f)$ is rather small, it seems that $g(\Xi^{train})$ does not contain enough points located near the borders of the set $g(\calX)$.
This is most probably due to the fact that the sampling strategy used to sample $\bfX$ is not a good sampling strategy to sample $g(\bfX)$, and maybe a re-sampling with a more suitable strategy would help to circumvent the problem.

Secondly for small $N^{train}$, we observe that GLI sometimes fails to converge towards the minimum of $\widehat{\calJ}$.
When it happens, we may not learn good functions $f$, that leads to large approximation errors $\widehat{e}_g(f)$ and $e_g(f)$.

\begin{figure}[H]
    \centering
    \includegraphics[page=2, width=0.99\textwidth]{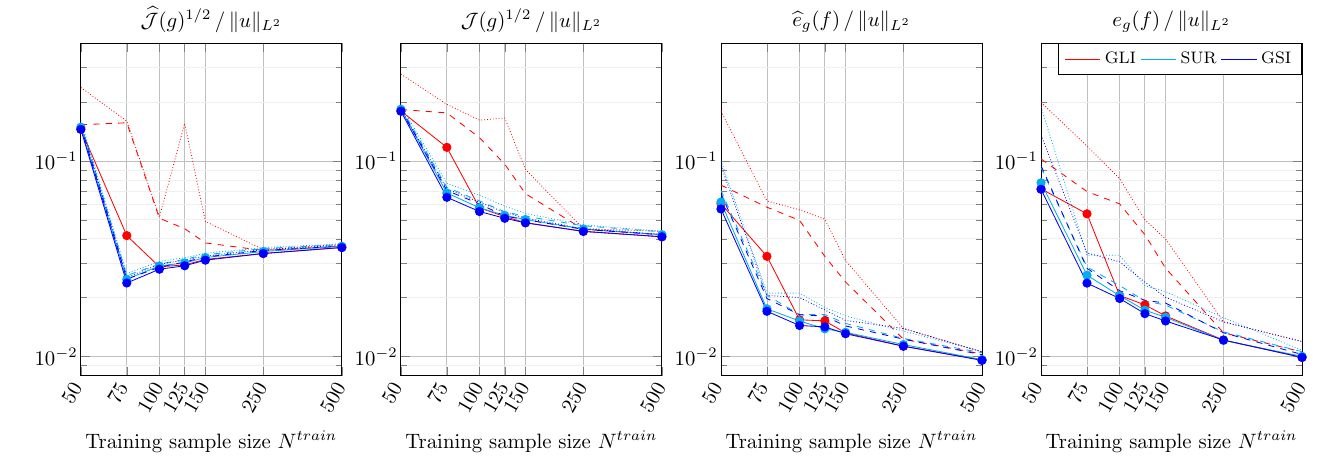}
    \caption[caption]{\footnotesize 
    Evolution of quantiles with respect to the size of the training sample for $u_3$ with $m=1$.
    The quantiles $50\%$, $90\%$ and $100\%$ are represented respectively by the continuous, dashed and dotted lines.
    }
    \label{fig:u3 vs ntrain}
\end{figure}

Let us continue with $u_3$, whose associated results are displayed in \Cref{fig:u3 vs ntrain}.
Firstly for all $N^{train}$, we observe that SUR and GSI perform at least as good as GLI, regarding all 4 monitored quantities.
Also, SUR and GSI yield very similar results, which can be interpreted as the minimizer of $\widehat{\calL}_1$ being near a local minimum of $\widehat{\calJ}$, if not near a global one.

Secondly for rather small $N^{train}$, we observe that SUR and GSI perform better than GLI on all $4$ quantities monitored.
We also observe that the results obtained with SUR and GSI are much less varying than with GLI over the $20$ realizations, especially regarding $\widehat{\calJ}(g)$ and $\calJ(g)$.

\begin{figure}[H]
    \centering
    \includegraphics[page=3, width=0.99\textwidth]{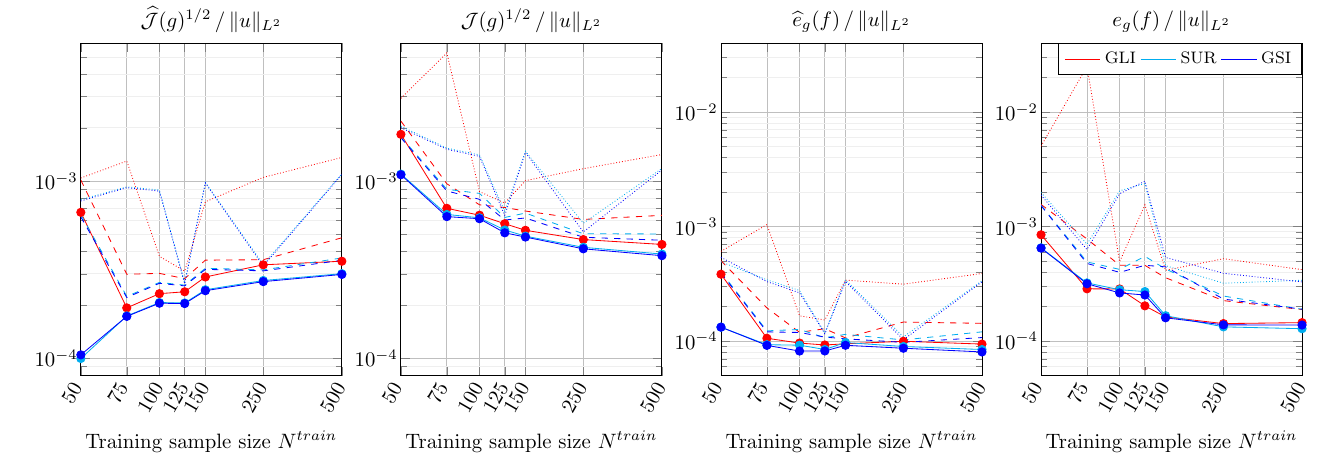}
    \caption[caption]{\footnotesize 
    Evolution of quantiles with respect to the size of the training sample for $u_4$ with $m=1$.
    The quantiles $50\%$, $90\%$ and $100\%$ are represented respectively by the continuous, dashed and dotted lines.
    }
    \label{fig:u4 vs ntrain}
\end{figure}

Let us finish with $u_4$, whose associated results are displayed in \Cref{fig:u4 vs ntrain}.
Firstly for all $N^{train}$, we again observe that SUR and GSI yield very similar results.
We also observe on the $100\%$ quantile that all 3 methods suffer from a rather significantly varying performance over the $20$ realizations, and that no method is always better than the others.
For example for $N^{train}=100$, the worst case errors for GLI are significantly better than with SUR and GSI, and vice versa for $N^{train}=75$.
However for $N^{train} \geq 75$, we observe on the $50\%$ and $90\%$ quantiles that GLI, SUR and GSI perform similarly.

Secondly for $N^{train}=50$, we observe that $\widehat{\calJ}(g)$ is often much lower for SUR and GSI than for GLI, while there is not much of a gap on $\calJ(g)$.
This means that, although SUR and GSI performed much better than GLI when it comes to minimizing $\widehat{\calJ}$, the minimizer of the latter is rather far from the one of $\calJ$.
Such behavior is often referred to as overfitting.

\subsection{The case of multiple features}
\label{subsec:results multiple features}

In this section we investigate the performances of the three learning procedures GLI, SUR and GLS when learning $m=2$ features to approximate $u_2$, $u_3$ and $u_4$.

\begin{figure}[H]
    \centering
    \includegraphics[page=4, width=0.99\textwidth]{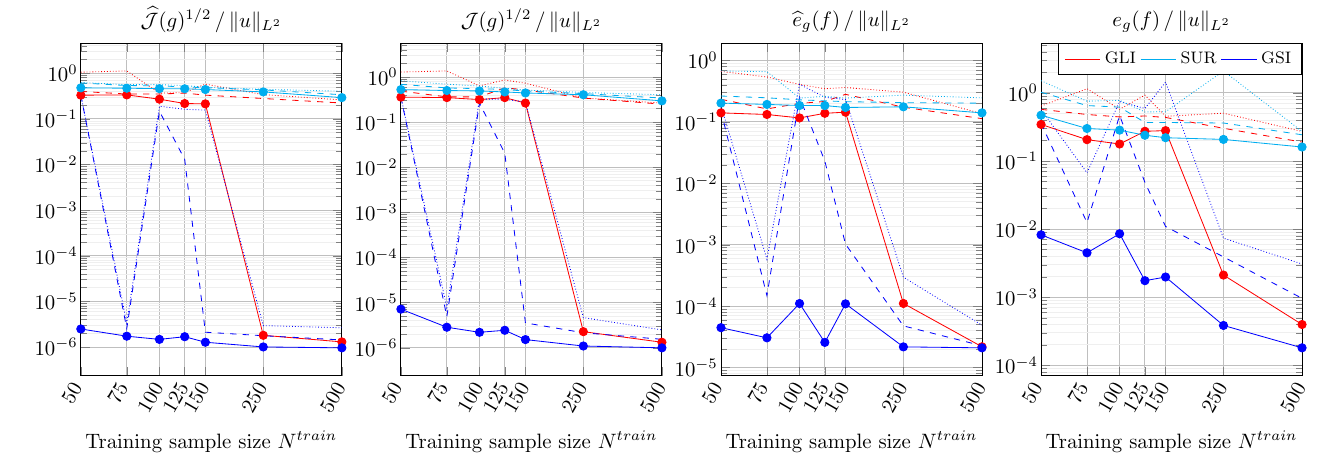}
    \caption[caption]{\footnotesize 
    Evolution of quantiles with respect to the size of the training sample for $u_2$ with $m=2$.
    The quantiles $50\%$, $90\%$ and $100\%$ are represented respectively by the continuous, dashed and dotted lines.
    }
    \label{fig:u2 vs ntrain m2}
\end{figure}

Let us start with $u_2$ for which results are displayed in \Cref{fig:u2 vs ntrain m2}.
Recall that this benchmark can be exactly expressed as a bivariate function of $m=2$ polynomial features of total degree $2$.

Firstly for all $N^{train}$, we observe that GSI outperforms GLI in all aspects, although the gap shrinks for large $N^{train}$.
In particular, it often yields the minimizer of $\widehat{\calJ}$,  which turns out to be $0$ as for $\calJ$.
We also observe that SUR and GSI now yield rather different results, which could be expected since the suboptimality results of our surrogates are much weaker for $m>1$ than for $m=1$.
In particular, SUR performs poorly here, although it appears that it yields a good initialization point for GSI.

Secondly for small $N^{train}$, we observe that GLI often fails to yield the minimizer of $\widehat{\calJ}$ while GSI is often successful, although there were realizations where GSI also failed. 

Thirdly for large $N^{train}$, we observe that GLI still sometimes fails to yield the minimizer of $\widehat{\calJ}$ while GSI is always successful.

\begin{figure}[H]
    \centering
    \includegraphics[page=5, width=0.99\textwidth]{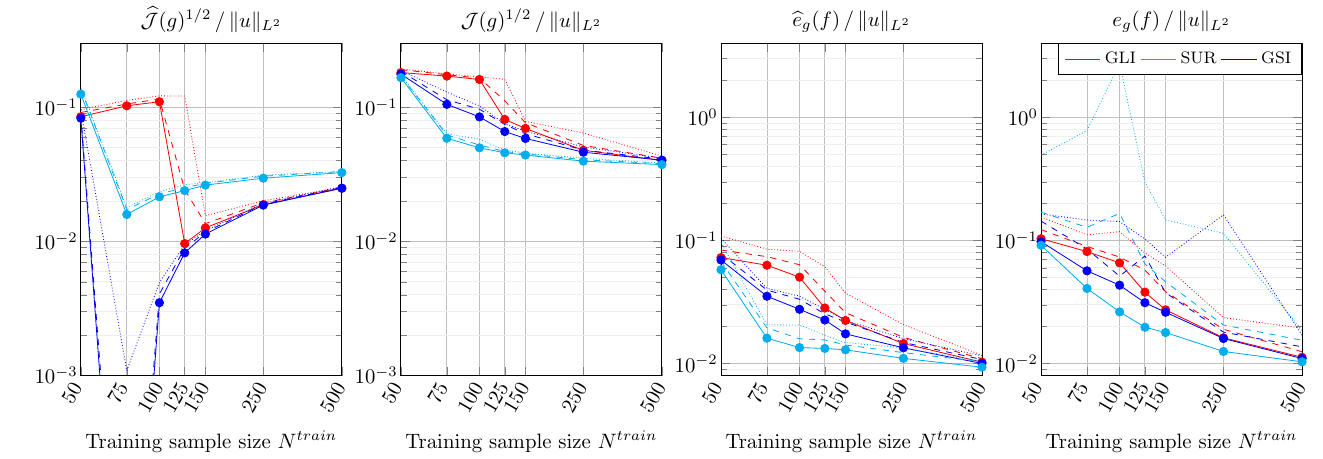}
    \caption[caption]{\footnotesize 
    Evolution of quantiles with respect to the size of the training sample for $u_3$ with $m=2$.
    The quantiles $50\%$, $90\%$ and $100\%$ are represented respectively by the continuous, dashed and dotted lines.
    }
    \label{fig:u3 vs ntrain m2}
\end{figure}

\begin{figure}[H]
    \centering
    \includegraphics[page=6, width=0.99\textwidth]{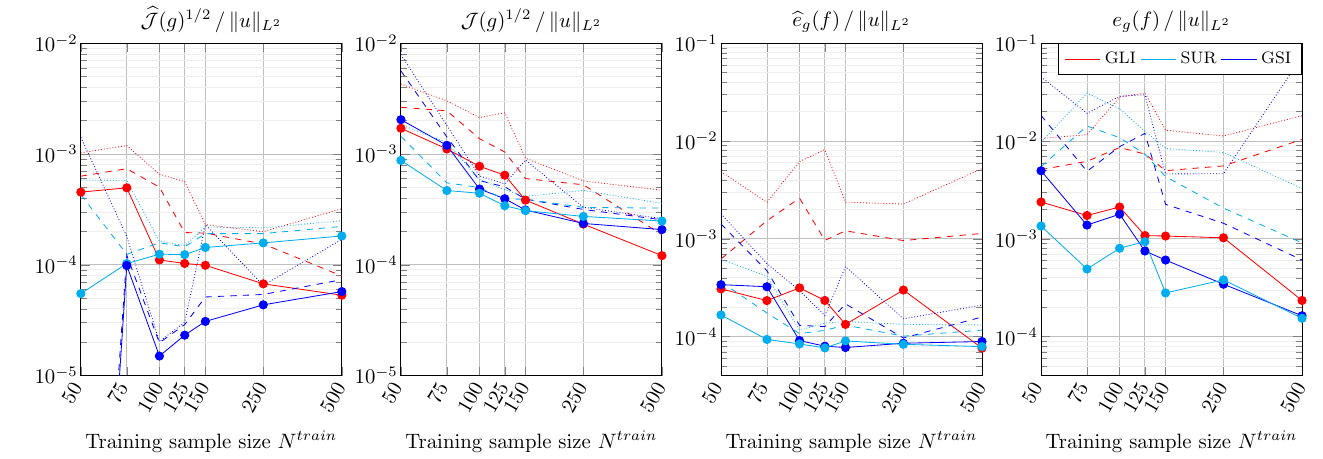}
    \caption[caption]{\footnotesize 
    Evolution of quantiles with respect to the size of the training sample for $u_4$ with $m=2$.
    The quantiles $50\%$, $90\%$ and $100\%$ are represented respectively by the continuous, dashed and dotted lines.
    }
    \label{fig:u4 vs ntrain m2}
\end{figure}

Let us continue with $u_3$ and $u_4$, whose associated results are displayed respectively in \Cref{fig:u3 vs ntrain m2,fig:u3 vs ntrain m2}.

Firstly we observe when comparing \Cref{fig:u3 vs ntrain,fig:u4 vs ntrain m2} that adding a feature did not improve much, if not deteriorated, the performances regarding $\calJ(g)$, $\widehat{e}_g(f)$ and $e_g(f)$, for all three procedures.

Secondly we observe that for rather small values of $N^{train}$, GSI performs much better at minimizing $\widehat{\calJ}(g)$ than SUR and GLI.
In particular for $u_3$ with $N^{train}=75$ and for $u_4$ with $N^{train}=50$, GSI mostly yields machine zero values for $\widehat{\calJ}(g)$.
However, the corresponding values for $\calJ(g)$ are much larger, which again translates same overfitting behavior.
This results in rather large approximation errors $\widehat{e}_g(f)$ and $e_g(f)$.
Hence for those benchmarks with $m=2$, it appears that the approaches SUR and GSI we proposed did not provide clear performance gain.

\subsection{General observations}
\label{subsec:results conclusion}

From what we observed in \Cref{subsec:results one features,subsec:results multiple features}, we can draw $3$ general observations.

Firstly, GSI mostly performed better than GLI for minimizing $\widehat{\calJ}(g)$, which is the Monte-Carlo estimation of $\calJ$ using training samples.
This is especially true for rather (but not too) small training sample sizes $N^{train}$.
For larger $N^{train}$ however, both GSI and GLI performed similarly.
Now the problem is that when $N^{train}$ is too small, the minimizer of $\widehat{\calJ}$ may be rather far from the one of $\calJ$, which is usually denoted as an overfitting behavior.
In such cases, the performance of GSI is similar to GLI, or even deteriorated.
We especially observed this overfitting behavior for $m=2$ on benchmarks $u_3$ and $u_4$ where the minimum of $\calJ$ was rather large.

Secondly, for $m=1$ feature, SUR yielded features that seemed to be near local minima of $\widehat{\calJ}$, and GSI did not yield much smaller values.
In particular when the minimum of $\widehat{\calJ}$ was $0$, as for $u_1$ and $u_2$, SUR always found it, while there were realizations where GLI did not.
In general, our observations suggest that our approach SUR, or its variant GSI, mostly learn significantly better features for (not too) small training sets.

Thirdly, for $m=2$ features, SUR yielded features that seemed rather far from local minima of $\widehat{\calJ}$, while GSI yielded significantly lower values.
When the problem was reducible to $m=2$ dimensions, GSI often yielded the correct dimension reduction, contrarily to GLI.
However, for benchmarks with not enough potential dimensionality reduction, it resulted in a rather strong overfitting behavior, resulting in mitigated results on $\calJ(f)$ and $e_g(f)$.

\section{Conclusion and perspectives}
\label{sec:conclusion}



In this work we investigated a nonlinear feature learning method for approximating high dimensional functions.
We considered a recently developed gradient-based method, which leverages Poincar\'e inequalities on nonlinear manifolds to derive an objective function.
One of the problem is that this objective function is a non-convex function of the feature map, which is problematic from an optimization point of view.
In order to circumvent this problem we introduced new surrogate quantities, which are quadratic with respect to the feature map, thus more suitable to optimization.
In particular, when taking feature maps in some linear space of nonlinear functions, minimizing these surrogates is equivalent to solving generalized matrix eigenvalue problems.

We first focused on the case of learning one feature, and we analyzed the sub-optimality of our surrogate.
We showed that our surrogate yields both upper and lower bounds on the objective function, for feature maps satisfying some deviation inequalities, typically polynomials.
We also extended this approach to the case of multiple features, for which we showed theoretical results similar to the case of one feature.

We finally illustrated our surrogate based methods on several numerical experiments.
We observed that our approach is especially efficient when the function to approximate admits a representation with a low-dimensional feature map, even for small sample sizes.
Using the feature map obtained by minimizing our surrogates as an initialization point of an iterative algorithm to minimize the Poincar\'e inequality based objective function improved the overall performance.


Let us mention three main perspectives to the current work.
The first perspective is to assess the performance of our approach on numerical examples with higher input dimension.

The second perspective is to investigate whether the surrogate we developed can effectively improve the optimization procedure for other classes of feature maps, such as diffeomorphisms-based feature maps, for which the theoretical guarantees of the surrogate can be much sharper than for polynomials, or such as feature maps based on low-rank tensor formats.

The third perspective is to extend our approach to the Bayesian inverse problem setting.
Indeed, recent works \cite{zahmCertifiedDimensionReduction2022,liPrincipalFeatureDetection2024} leveraged gradient-based functional inequalities to derive nonlinear dimension reduction methods for approximating the posterior distribution.
Extending our surrogate methods may improve the learning procedure of nonlinear feature in such setting.

\begin{appendices}

\section{Proofs for \Cref{sec:Large and small deviations inequalities}}
\label{apx:proof deviation ineq}

\subsection{Proof of \Cref{prop:large deviation remez uniform}}
\label{subsec:proof prop:large deviation remez uniform}

Let $h$ satisfying \eqref{equ:remez inequality multivariate} with constants $k,A$.
If $h=0$ then the desired inequality holds as $q_h=0$ and $\Proba{|h(\bfX)|>0} = 0$. 
Otherwise, using the fact that $q_{|h|^{1/k}}^k=q_h$ and that $\Proba{|h(\bfX)| > q_h t} \leq \Proba{|h(\bfX)| \geq q_h t}$, we have from \cite{fradeliziConcentrationInequalities$s$concave2009} that for all $t>A^k$,
\begin{equation}
\label{equ:proof large deviation}
  \Proba{|h(\bfX)| > q_h t}
  \leq
  \left\{
  \begin{aligned}
    (1-(1-2^{-s})A^{-1}t^{1/k})^{1/s}_+,
    \quad &s \in (0,1/d], 
    \\
    \exp(-\log(2) A^{-1} t^{1/k}),
    \quad &s=0,
    \\
    (2^{-s}-1)^{1/s} A^{-1/s} t^{1/sk},
    \quad &s<0.
  \end{aligned}
  \right.
\end{equation}
Now, we will  show that the previous equation actually holds for all $t>1$, then we will weaken the bound so that it holds for all $t>0$.
Recall that by definition $k,A \geq 1$.

\paragraph{Case $0<s\leq 1/d$.}
Let $t \in (1, A^k]$.
Since $t > 1$ and $\Proba{|h(\bfX)| \leq q_h} \geq 1/2$, we have that
\[
  \Proba{|h(\bfX)| >  q_h t}
  = 1-\Proba{|h(\bfX)| \leq q_h t}
  \leq 1/2.
\]
Now since $x \mapsto (1-(1-2^{-s})A^{-1}x^{1/k})^{1/s}_+$ is decreasing and $t \leq  A^k$, we have that
\[
  (1-(1-2^{-s})A^{-1} t^{1/k})^{1/s}_+
  \geq
  (1-(1-2^{-s}))^{1/s}_+
  = 
  1/2,
\]
which yields that \eqref{equ:proof large deviation} holds for all $t> 1$.
Now let $t>0$ and $\eta = \overline{\eta}_{A,s} := A(1-2^{-s})^{-1}$.
We have that $(1-t^{1/k} / \eta ) \leq (1- (t^{1/k} - 1) / \eta )$, and if $t\in(0,1]$ we have that $1 \leq (1-(t^{1/k} - 1) / \eta )$.
Thus, we have the desired result, for all $t>0$,
\[
  \Proba{|h(\bfX)| > \overline{q}_h t}
  \leq  (1-(t^{1/k} - 1) / \eta ).
\]

\paragraph{Case $s=0$.}
Let $t\in (1, A^k]$.
As for the case $s>0$, we have $\Proba{|h(\bfX)| > q_ht} \leq 1/2$.
Now since $x \mapsto \exp(-\log(2)A^{-t} x^{1/k})$ is decreasing and $t \leq A^k$, we have that 
\[
  \exp(-\log(2)A^{-t} t^{1/k}) \geq  \exp(-\log(2)) = 1/2,
\]
which yields that \eqref{equ:proof large deviation} holds for all $t> 1$.
Now, let $t>0$ and $\eta = \overline{\eta}_{A, s} :=A / \log(2)$.
We have that $\exp(- t^{1/k}/\eta) \leq \exp(-(t^{1/k}- 1)/\eta)$ and if $t\in(0,1]$ we have that $1 \leq \exp(-(t^{1/k}- 1)/\eta)$.
Thus, we have the desired result for, for all $t>0$,
\[
  \Proba{|h(\bfX)| > q_ht}
  \leq  \exp(-(t^{1/k}- 1)/\eta).
\]

\paragraph{Case $s<0$.}
Let $t\in (1, A^k]$.
As for the case $s>0$, we have $\Proba{|h(\bfX)| \geq q_h t} \leq 1/2$.
Now since $s<0$, the function $x \mapsto (2^{-s} - 1)^{1/s}A^{-1/s}x^{1/sk}$ is decreasing, and since $t \leq A^k$ and $s<0$ we have that $(A^{-1}t^{1/k})^{1/s} \geq 1$ and thus
\[
  (2^{-s} - 1)^{1/s}A^{-1/s}t^{1/sk} \geq  (2^{-s} - 1)^{1/s} \geq 1/2,
\]
which yields that \eqref{equ:proof large deviation} holds for all $t> 1$.
Now let $t>0$ and $\eta = \overline{\eta}_{A,s} := \max(1,A / (2^{-s}-1))^{-1/s}$.
We have that $(2^{-s}-1)^{1/s} A^{-1/s} t^{1/sk} \leq \eta t^{1/sk}$ and if $t\in(0,1]$ we have that $\eta t^{1/sk} \geq 1$ since $s<0$.
Thus, we have the desired result for, for all $t>0$,
\[
  \Proba{|h(\bfX)| > \overline{q}_h t}
  \leq  \eta t^{1/sk}.
\]

\subsection{Proof of \Cref{prop:bounded uniform quantiles}}
\label{subsec:proof prop:bounded uniform quantiles}

Let us first show the desired left inequality.
Let $\bfX$ a random variable on $\bbR^d$ and let $h$ satisfying \eqref{equ:remez inequality multivariate}.
Assume that $\Proba{h(\bfX) = 0} > 0$, then there exists $\bfx_0$ such that $h=0$ in a neighborhood of $\bfx_0$.
Then for any $\bfx \neq \bfx_0$, there exists $\varepsilon>0$ small enough such that $\sup_{I_{\varepsilon}} |h| = 0$ with $I_{\varepsilon} := \{t \bfx_0 + (1-t) \bfx: t\in[0,\varepsilon]\} \subset I_{1}$, thus from \eqref{equ:remez inequality multivariate} it holds 
\[
  \sup_{I_1}|h| 
  \leq 
  \Big(
    \frac{A \|\bfx - \bfx_0\|}{\varepsilon \|\bfx - \bfx_0\|}
  \Big)^k
  \sup_{I_{\varepsilon}}|h|
  = \Big(
    \frac{A }{\varepsilon}
    \Big)^k
  \sup_{I_{\varepsilon}}|h|
  = 0,
\]
thus $h(\bfx)=0$, which contradicts $h$ being not identically zero, thus $\Proba{h(\bfX) = 0} = 0$.
Then, applying the Markov's inequality on $|h(\bfX)|$, which is a positive real random variable from the previous statement, we obtain for all $\omega \in (0,1)$ and $p>0$,
\[
  Q_{|h(\bfX)|}(\omega) 
  \leq (1-\omega)^{-1/p} \Expe{|h(\bfX)|^p}^{1/p}
  = (1-\omega)^{-1/p} \|h\|_{L^p_{\mu}}.
\]
Finally, taking $\omega=1/2$ yields the desired left inequality.
Let us now show the desired right inequality.
Let us denote $q = Q_{|h(\bfX)|}(1/2)$, which satisfies $q = Q_{|h(\bfX)|^{1/k}}(1/2)^k$.

\paragraph{Case $s\in (0,1/d]$.}

From \cite[Corollary 10]{fradeliziConcentrationInequalities$s$concave2009}, since $0 < 1 - 2^{-s} < 1$ we have that 
\[
  \Proba{
    |h(\bfX)|^{1/k} \geq A q^{1/k} (1-2^{-s})^{-1}
  } \leq 0,
\]
hence $\|h^{1/k}\|_{L^{\infty}_{\mu}} \leq A q^{1/k} (1-2^{-s})^{-1}$.
Now since $\|h^{1/k}\|_{L^{\infty}_{\mu}} = \|h\|_{L^{\infty}_{\mu}}^{1/k} \geq \|h\|_{L^{p}_{\mu}}^{1/k}$ for all $p\in (0,+\infty]$, we obtain the desired inequality, for all $p\in(0,+\infty]$,
\[ 
  q \geq \|h\|_{L^{p}_{\mu}} A^{-k} (1-2^{-s})^{k}.
\]

\paragraph{Case $s\geq 0$.}
From \cite[Corollary 10]{fradeliziConcentrationInequalities$s$concave2009} we have that for all $p>0$,
\[
  \||h|^{1/k}\|_{L^p_{\mu}}
  \leq A q^{1/k} (1 + 2^{p} \Gamma(p+1))^{1/p},
\]
where $\Gamma$ denotes the classical Gamma function.
Then for any $p\in (0,+\infty)$, by applying the above result for $\tilde p = pk>0$ and noting that $\||h|^{1/k}\|_{L^{pk}_{\mu}} = \|h\|_{L^{p}_{\mu}}^{1/k}$, we obtain for all $p\in (0,+\infty)$,
\[
  q \geq \|h\|_{L^{p}_{\mu}} A^{-k} (1 + 2^{pk} \Gamma(pk+1))^{-1/p}.
\] 
Now assume that $p\geq 1/k$.
Then applying \Cref{lem:bound gamma function} below with $y=pk\geq 1$ yields the desired result,
\[
  q \geq \|h\|_{L^{p}_{\mu}} A^{-k} (3pk)^{-k}.
\]

\begin{lemma}
\label{lem:bound gamma function}
  For all $y\geq 1$ it holds $(1+2^y \Gamma(y+1))^{1/y} \leq 3y$.
\end{lemma}
\begin{proof}
  Let $y\geq 1$.
  Since $2^y \Gamma(y+1) \geq 2$ we have 
  \[
    (1+2^y \Gamma(y+1))^{1/y} 
    \leq 2 \Gamma(y+1)^{1/y} (1 + \frac{1}{2^y \Gamma(y+1)})^{1/y}
    \leq 3 \Gamma(y+1)^{1/y}.
  \]
  Then, \cite[Theorem 1.6]{batirInequalitiesGammaFunction2009} states that $\Gamma(y+1)<y^y e^{-y} \sqrt{2\pi(y+b)}$ with $b:=e^2 / 2\pi -1 >0$, which yields
  \[
    (1+2^y \Gamma(y+1))^{1/y} 
    \leq 3 y e^{-1} (2\pi(y+b))^{1/2y}
    = 3 y e^{-1} h(y)^{1/2}.
  \]
  with $h:x \mapsto (2\pi(x+b))^{1/x}$.
  Now for $x\geq 1$, after calculation we obtain
  \[
    h'(x)
    = \frac{h(x)}{x^2(x+b)}
    \Big(
      x - (x+b) \log(2\pi(x+b))
    \Big),
  \]
  and since $2\pi(x+b) = 2\pi(x-1) + e^2 \geq e$, we have $\log(2\pi(x+b)) > 1$ thus $h'(x)<0$ and $h(x) \leq h(1)=e^2$, which yields the desired result.
\end{proof}

\paragraph{Case $s<0$.}
We will use \cite[Corollary 11]{fradeliziConcentrationInequalities$s$concave2009}.
However, note that it contains a sign error.
After correcting the sign error, we have that for all $p\in(0,-1/s)$,
\[
  \||h|^{1/k}\|_{L^p_{\mu}}
  \leq A q^{1/k} 
  \Big(
    1 - (2^{-s}-1)^{1/s} \frac{p}{p + 1/s} 
  \Big)^{1/p}.
\]
Then for any $p\in (0,-1/sk)$, by applying the above result for $\tilde p = pk \in (0, -1/s)$ and noting that $\||h|^{1/k}\|_{L^{pk}_{\mu}} = \|h\|_{L^{p}_{\mu}}^{1/k}$, we obtain the desired inequality, for all $p\in (0,-1/sk)$,
\[
  q \geq 
  \|h\|_{L^{p}_{\mu}} A^{-k} 
  \Big(1 - \frac{(2^{-s}-1)^{1/s}}{1 + 1/spk} \Big)^{-1/p}.
\]

\section{Proofs for \Cref{sec:one feature}}
\label{apx:proof one feature}

\subsection{Proof of \Cref{prop:bound J by L}}
\label{subsec:proof prop:bound J by L}

Let $g\in\calG_1$ and define for all $\alpha >0$ the event $E(\alpha) := (\|\nabla g(\bfX)\|_2^2 < \alpha)$.
First using $1 = \mathbbm{1}_{E(\alpha)} + \mathbbm{1}_{\overline{E(\alpha)}}$ and $\|\Piperp_{\nabla g(\bfX)} \nabla u(\bfX)\|_2^2 \leq  \|\nabla u(\bfX)\|^2_2$, we obtain that for all $\alpha >0$,
\[
  \calJ(g) 
  \leq \Expe{\|\Piperp_{\nabla g(\bfX)} \nabla u(\bfX)\|_2^2 \mathbbm{1}_{\overline{E(\alpha)}}}
  + \Expe{\|\nabla u(\bfX)\|_2^2 \mathbbm{1}_{E(\alpha)}}.
\]
To upper bound the second term in the right-hand side of the above inequality, we use Holder's inequality and the fact that $\Expe{\|\nabla u(\bfX)\|_2^{2p_u}} \leq 1$.
To upper bound the first term, we first use \Cref{lem:norm projection} under the event $E(\alpha)$, then we use $\mathbbm{1}_{E(\alpha)} \leq 1$, and finally we recognize the definition of $\calL_1$ from \eqref{equ:def L one feature}.
As a result, we obtain for all $\alpha >0$,
\[
  \calJ(g)
  \leq \alpha^{-1} \calL_1(g) 
  + \Proba{\|\nabla g(\bfX)\|_2^2 < \alpha}^{1/p}.
\]
Now, from \Cref{assump:bounded uniform A k} we can apply \Cref{prop:small deviation remez uniform} on $h : \bfx \mapsto \|\nabla g(\bfx)\|_2^2$.
Moreover from \Cref{assump:full rank jacobian as}, $h$ is not the zero function, thus we have from \Cref{prop:bounded uniform quantiles} that $q_h>0$.
We then to obtain for all $\alpha >0$,
\[
  \Proba{\|\nabla g(\bfX)\|_2^2 < \alpha} 
  \leq \underline{\eta}_{A, s} (\alpha / q_h)^{1/k}.
\]
Now, denoting $\xi := (\underline{\eta}_{A,s} q_h^{-1/k})^{1/p}$, we obtain for all $\alpha >0$,
\[
  \calJ(g) \leq 
  \xi \big( \xi^{-1}\calL_1(g) \alpha^{-1} + \alpha^{1/p k}\big).
\]
We can now study the function $\alpha \mapsto a\alpha^{-1} + \alpha^b$ with $a=\xi^{-1}\calL_1(g)$ and $b=1/pk$ to optimize the above upper bound.
This is done in \Cref{lem:technical lem 1} below, and yields
\[
  \calJ(g) \leq 
  2 \xi (\xi^{-1}\calL_1(g))^{\frac{1}{1+pk}}
  = 2\xi^{\frac{pk}{1+pk}} \calL_1(g)^{\frac{1}{1+pk}}
  = 2 \underline{\eta}_{A, s}^{\frac{k}{1+pk}} 
  q_h^{-\frac{1}{1+pk}} \calL_1(g)^{\frac{1}{1+pk}}.
\]
Now, using \Cref{prop:bounded uniform quantiles} for any $p_1 \geq 1$ such that $-1 < skp_1 \leq +\infty$, we can upper bound $q_h^{-1}$ using $\Expe{\|\nabla g(\bfX)\|_2^{2p_1}}^{-1/p_1} \leq \underline{\nu}_{\calK_1, p_1}^{-1}$, depending on the value of $s$.
Firstly if $s\in (0,1/d]$ then
\[
  \calJ(g) \leq 
  2(\underline{\eta}_{A,s} \overline{\eta}_{A,s})^{\frac{k}{1+pk}}
  \underline{\nu}_{\calK_1,p_1}^{-\frac{1}{1 + pk}}
  \calL_1(g)^{\frac{1}{1 + pk}}.
\]
Secondly if $s \geq 0$ then
\[
  \calJ(g) \leq 
  2(3\underline{\eta}_{A,s} A k p_1
  )^{\frac{k}{1+pk}}
  \underline{\nu}_{\calK_1,p_1}^{-\frac{1}{1 + pk}}
  \calL_1(g)^{\frac{1}{1 + pk}}.
\]
Thirdly if $s \in (-1/k, 0)$ then
\[
  \calJ(g) \leq 
  2(
  \underline{\eta}_{A,s}
  A (1- \frac{(2^{-s}-1)^{1/s}}{1+ 1/skp_1} )^{\frac{1}{kp_1}}
  )^{\frac{k}{1+pk}}
  \underline{\nu}_{\calK_1,p_1}^{-\frac{1}{1 + pk}}
  \calL_1(g)^{\frac{1}{1 + pk}}.
\]
Combining the three previous inequalities and the fact that $s$-concave measures with $s\in(0,1/d)$ are also $0$-concave measures yields the desired result.

\begin{lemma}
  \label{lem:technical lem 1}
  Let $a\geq 0$ and $b>0$. Then 
  \begin{equation}
    a^{\frac{b}{1 + b}}
    \leq \inf_{x>0} (a x^{-1} + x^b) 
    \leq 2 a^{\frac{b}{1 + b}}.
  \end{equation}
\end{lemma}
\begin{proof}
  If $a=0$ then the result is straightforward. 
  Now assume that $a>0$ and let 
  \[
    h : x \mapsto a x^{-1} + x^b.
  \]
  We have $h\in \calC^{\infty}(\bbR^*_+)$ and $h(x) \rightarrow +\infty$ as $x\rightarrow \{0, + \infty\}$, thus $h$ admits at least 1 global minimum.
  Now, we have $h'(x) = - a x^{-2} + b x^{b-1}$, thus the global minimum of $h$ is at $x = (\frac{a}{b})^{\frac{1}{1+b}}$ and has value $(1+b) (\frac{a}{b})^{\frac{b}{1+b}}$.
  In other words, we can write
  \[
    \inf_{x>0} h(x) = v(b)a^{\frac{b}{1+b}},
  \]
  where $v(b) = b\mapsto (1+b) b^{-\frac{b}{1 + b}}$. 
  Now we have $v\in \calC^{\infty}(\bbR^*_+)$ and $v(x) \rightarrow 1$ as $x\rightarrow \{0, + \infty\}$.
  Moreover,
  \begin{equation}
    v'(x) =
    -\frac{\log(x)}{1+x} x^{-\frac{x}{1+x}},
  \end{equation}
  thus $v$ is strictly increasing on $(0,1)$ and strictly decreasing on $(1,+\infty)$. Thus it holds $1 = v(0) \leq v(b) \leq v(1) = 2$, which yields the desired result.
\end{proof}

\subsection{Proof of \Cref{prop:bound L by J}}
\label{subsec:proof prop:bound L by J}

Let $g\in\calG_1$ and let us first assume that $\Expe{\|\nabla g(\bfX)\|_2^{2p_1}} = 1$.
Define for all $\alpha >0$ the event $E(\alpha) := (\|\nabla g(\bfX)\|_2^2 > \alpha)$.
First using $1 = \mathbbm{1}_{E(\alpha)} + \mathbbm{1}_{\overline{E(\alpha)}}$ and $\|\Piperp_{\nabla u(\bfX)} \nabla g(\bfX)\|_2^2 \leq  \|\nabla g(\bfX)\|^2_2$, we obtain that for all $\alpha >0$,
\[
  \calL_1(g)
  \leq \Expe{
    \|\nabla u(\bfX)\|_2^2 
    \|\Piperp_{\nabla u(\bfX)} \nabla g(\bfX)\|_2^2
    \mathbbm{1}_{\overline{E(\alpha)}}
  }
  + \Expe{
    \|\nabla u(\bfX)\|_2^2 
    \|\nabla g(\bfX)\|_2^2
    \mathbbm{1}_{E(\alpha)}
  }.
\]
To upper bound the second term in the right-hand side of the above inequality, we use Holder's inequality with the facts that $\Expe{\|\nabla g(\bfX)\|_2^{2p_1}} = 1 \geq  \Expe{\|\nabla u(\bfX)\|_2^{2p_u}}$ and $r^{-1} = 1 - p_u^{-1} - p_1^{-1}>0$ by assumption.
To upper bound the first term, we first use \Cref{lem:norm projection} under the event $E(\alpha)$, then we use $\mathbbm{1}_{E(\alpha)} \leq 1$, and finally we recognize the definition of $\calJ(g)$ from \eqref{equ:def of J}.
As a result, we obtain for all $\alpha>0$,
\[
  \calL_1(g) \leq
  \alpha \calJ(g) + \Proba{\|\nabla g(\bfX)\|_2^2 > \alpha}^{1/r}.
\]
Now from \Cref{assump:bounded uniform A k} we can apply \Cref{prop:small deviation remez uniform} on $h : \bfx \mapsto \|\nabla g(\bfx)\|_2^2$.
Moreover from \Cref{assump:full rank jacobian as}, $h$ is not the zero function thus \Cref{prop:bounded uniform quantiles} yields $0 < q_h \leq 2 \|h\|_{L^{p_1}_{\mu}} = 2$.

Firstly if $0<s<1/d$, which implies that $\bfX$ has compact support, then from \Cref{prop:bounded uniform quantiles} we have $\|h\|_{L^{\infty}_{\mu}} \leq 2 \overline{\eta}_{A,s}^k$, thus taking $\alpha =2 \overline{\eta}_{A,s}^k$ yields $\Proba{\|\nabla g(\bfX)\|_2^2 > \alpha} = 0$.
Hence we obtain the desired result
\[
  \calL_1(g)\leq \gamma_2 \calJ(g),
  \quad
  \gamma_2 := 2\overline{\eta}_{A,s}^{k}.
\]
Secondly if $s=0$ then by \Cref{prop:large deviation remez uniform} and using $q_h\leq 2$ we have for all $\alpha >0$,
\[
  \Proba{\|\nabla g(\bfX)\|_2^2 > \alpha}
  \leq 
  \exp(-((\alpha / 2)^{1/k} - 1)/\overline{\eta}_{A,s}) 
  = \exp(1/\overline{\eta}_{A,s}) \exp(-r\beta)
\] 
where $\beta=(\overline{\eta}_{A,s} r)^{-1} (\alpha / 2)^{1/k}>0$.
Then since $A,r\geq 1$ we have $r \overline{\eta}_{A,s} = rA / \log(2) \geq 1$ and $\exp(1/\overline{\eta}_{A,s} r) = 2^{1/Ar} \leq 2$, thus we have that $\exp(1/\overline{\eta}_{A,s} r) \leq 2(\overline{\eta}_{A,s}r)^k$.
Hence for all $\beta >0$,
\[
  \calL_1(g)\leq
  2 (\overline{\eta}_{A,s} r)^{k} \calJ(g) \beta^{k}
  + \exp(1/\overline{\eta}_{A,s} r) \exp(-\beta)
  \leq \kappa (a \beta^k + \exp(-\beta)),
\]
with $\kappa := 2 (\overline{\eta}_{A,s} r)^{k}$ and $a:= \calJ(g)$.
Now, since $\Expe{\|\nabla u(\bfX)\|_2^{2p_u}} \leq 1$ with $p_u > 1$ we have $a\leq 1$.
Then one can show that the value of $\beta$ which minimizes the right-hand side behaves like $\log(a^{-1})$ as $a\rightarrow 0$.
Hence, we evaluate the right-hand side at $\beta = \log(a^{-1}) \geq 0$ to obtain $\calL_1(g)\leq \kappa a (1 + \log(a^{-1})^k)$.
In other words, we have the desired result,
\[
  \calL_1(g)
  \leq \gamma_2
  \calJ(g) (1 + |\log(\calJ(g))|^k),
  \quad
  \gamma_2 := 2 (\overline{\eta}_{A,s} r)^{k}.
\] 
Thirdly if $s\in(-1/k, 0)$ then by \Cref{prop:large deviation remez uniform} and using $q_h\leq 2^{1/p_1}$ we have for all $\alpha >0$,
\[
\Proba{\|\nabla g(\bfX)\|_2^2 > \alpha} 
\leq \overline{\eta}_{A,s} (\alpha/2^{1/p_1})^{1/sk}
= \kappa^{r} \beta^{-r}
\]
with $\beta := \alpha^{-1/srk}$ and $\kappa =\overline{\eta}_{A,s}^{1/r}2^{-1/srkp_1}$, with $\kappa < (2\overline{\eta}_{A,s})^{1/r}$ since $p_1^{-1}> -sk$.
We then have that for all $\beta >0$,
\[
  \calL_1(g)\leq
  \calJ(g) \beta^{-srk} + \kappa \beta^{-1}
  = \calJ(g) (a \beta^{-1} + \beta^{b})
\]
with $a=\kappa \calJ(g)^{-1} \geq 0$ and $b = -srk >0$.
Hence applying \Cref{lem:technical lem 1} and using $\frac{-sk}{1-srk} \leq r^{-1} \leq 1$ we obtain the desired result,
\[
  \calL_1(g)\leq
  \calJ(g) 2 a^{\frac{b}{1+b}} 
  = 2(2\overline{\eta}_{A,s})^{\frac{-sk}{1-srk}} 
  \calJ(g)^{\frac{1}{1-srk}}
  \leq \gamma_2  \calJ(g)^{\frac{1}{1-srk}},
  \quad
  \gamma_2 := 4\overline{\eta}_{A,s}^{\frac{1}{r}}.
\]
Finally, if $\Expe{\|\nabla g(\bfX)\|_2^{2p_1}} \neq 1$, then we can apply the same reasoning to $\tilde g : \bfx \mapsto \alpha^{-1/2}_g g(\bfx)$ with $\alpha_g := \Expe{\|\nabla g(\bfX)\|_2^{2p_1}}^{1/p_1}$.
Then using $ \calJ(g)=\calJ(\tilde g)$, $\calL_1(g) = \alpha_g \calL_1(\tilde g)$ and $\alpha_g \leq \overline{\nu}_{\calK_1, p_1}$, we obtain the desired result.

\subsection{Proof of \Cref{cor:suboptimality uniform}}
\label{subsec:proof cor:suboptimality uniform}

Since $\bfX$ is uniformly distributed on a compact convex subset of $\bbR^d$, its law is $s$-concave with $s=1/d$.
Also, since $\calG_1$ contains only non-constant polynomials of total degree at most $\ell+1$, we have that $\bfx\mapsto \|\nabla g(\bfx)\|_2^2$ is a non zero polynomial of degree at most $2\ell$ for any $g\in\calG_1$.
As a result, $\calG_1$ satisfies \Cref{assump:bounded uniform A k} with constants $k=2\ell$ and $A=4$, as well as \Cref{assump:full rank jacobian as}.
Now, applying \Cref{prop:bound L by J,prop:bound J by L} with $p_u=+\infty$ since $\bfX$ has compact support and $p_1=1$, we obtain for all $g\in\calG_1$,
\[
  \gamma_2^{-1} \calL_1(g)
  \leq \calJ(g)
  \leq \gamma_1 \calL_1(g)^{\frac{1}{1+k}},
  \quad \gamma_1 = 2(\underline{\eta}_{A,s} A \min\{ 3k,\frac{1}{1-2^{-s}})\})^{\frac{k}{1+k}},
  \quad  \gamma_2 = 2\overline{\eta}_{A,s}^{k}.
\]
Moreover $\underline{\eta}_{A,s} = A (1-2^{-s})s^{-1} \leq 4\log(2)\leq 4$ and $(1-2^{-s})^{-1} \leq 2d$, thus $\gamma_1 \leq 2 (16\min\{3k, 2d\})^{k/(1+k)}$ and $\gamma_2 \leq 2(8d)^k$.
Now, using \Cref{prop:sub optimality} with $p_u=+\infty$ and $p_1=1$ yields that any minimizer $g^*$ of $\calL_1$ over $\calG_1$ satisfies
\[
  \calJ(g^*)
  \leq \gamma_3 \inf_{g\in \calG_1} \calJ(g)^{\frac{1}{1+k}},
  \quad
  \gamma_3 = \gamma_1 \gamma_2^{\frac{1}{1+k}}
  = 2^{1+\frac{1}{1+k}} (A^3 s^{-1} \min\{3k, \frac{1}{1-2^{-s}}\})^{\frac{k}{1+k}}.
\]
Finally, using $(1-2^{-s})^{-1} \leq 2d$ and $k=2\ell$ we obtain the desired inequality $\gamma_3 \leq 4(2^7 d \min\{3\ell, d\})^{\frac{2\ell}{1+2\ell}}$.

\subsection{Proof of \Cref{prop:sharpness of rate with poly degree}}
\label{subsec:proof prop:sharpness of rate with poly degree}

Let $\bfX \sim \calU((0,1)^2)$, $\ell \in \bbN^*$ and $u(\bfx) = x_1$.
Let $a\in[0,1]$ and $g^a$ as defined in \Cref{prop:sharpness of rate with poly degree}.
After calculation we obtain
\[
  \calL_1(g^a)
  = \Big(1 + \frac{\Expe{X_1^{2\ell}}}{a^{2\ell} \Expe{X_2^{2\ell}}} \Big)^{-1},
  \quad
  \calJ(g^a)
  = \Expe{\Big(1 + \frac{X_1^{2\ell}}{a^{2\ell} X_2^{2\ell}} \Big)^{-1}}.
\]
Using the fact that $X_1$ and $X_2$ are i.i.d we obtain the desired result on $\calL_1$.
Then, let us write $\calJ(g^a) = \Expe{\big(1 + (Y/a)^{2\ell}\big)^{-1}}$ with $Y := X_1 / X_2$.
We know that the ratio of two independent uniform random variables has density $\rho_{Y} : t\mapsto \frac{1}{2}\mathbbm{1}_{(0,1)}(t) + \frac{1}{2t^2}\mathbbm{1}_{(1, +\infty)}(t)$.
Then, using the fact that $\rho_{Y}(t) \geq \frac{1}{2}\mathbbm{1}_{(0,1)}(t)$, 
applying a change of variable 
and using the fact that $\int_0^{1/a} \frac{1}{1 + t^{2\ell}} dt \geq \frac{1}{2}$ for $a$ small enough, since $\int_0^b \frac{1}{1 + t^{2\ell}} dt \xrightarrow[b\rightarrow +\infty]{} \frac{\pi / 2\ell}{\sin(\pi / 2\ell)} \geq 1$, we obtain
\[
\calJ(g^a) 
\geq  
\frac{1}{2} \int_0^1 \frac{1}{1 + (t/a)^{2\ell}} dt
= \frac{a}{2} \int_0^{1/a} \frac{1}{1 + t^{2\ell}} dt
\geq \frac{a}{4},
\]
which is the desired result on $\calJ(g^a)$.
As a result for $a\rightarrow 0$, since $\calL_1(g^a) \sim a^{2\ell}$, we have that $\calL_1(g^a)^{1/2\ell} \lesssim \calJ(g^a)$.
Finally, applying \Cref{prop:bound J by L} with $s=1/d$ and $p=1$ yields that $\calJ(g^a) \lesssim \calL_1(g^a)^{1 / (1+2\ell)}$.
These two last inequalities yield the last desired result.

\section{Proofs for \Cref{sec:multiple features}}
\label{apx:proof multiple features}

\subsection{Proof of \Cref{lem:bounds on norm w}}
\label{subsec:proof lem:bounds on norm w}

Let $1\leq j\leq m$, $g\in \calC^1(\calX, \bbR^m)$ and $\bfx\in\calX$.
Firstly, the definition of $w_{g,j}$ in \eqref{equ:def w, v} directly yields the first desired inequality $\|w_{g,j}(\bfx)\|_2^2 \leq \|\nabla g_j(\bfx)\|_2^2$.
Secondly, assume that $\sigma_{m}(\nabla g(\bfx)) > 0$, otherwise the second desired inequality holds trivially.
Let us first assume that $j=m$, let $M := \nabla g(\bfx)^T \nabla g(\bfx) \in \bbR^{m\times m}$ with corresponding blocks $M_1 := \nabla g_{-j}(\bfx)^T \nabla g_{-j}(\bfx) \in \bbR^{(m-1)\times (m-1)}$, $M_2 := \nabla g_{-j}(\bfx)^T \nabla g_j(\bfx) \in \bbR^{(m-1) \times 1}$ and $M_{3} := \nabla g_j(\bfx)^T \nabla g_j(\bfx) \in \bbR$ such that
\[
M =
\begin{pmatrix}
  M_1 & M_2 \\
  M_2^T & M_3
\end{pmatrix}.
\]
Since $\sigma_m(\nabla g(\bfx))>0$ we have that $\sigma_{m-1}(\nabla g_{-j}(\bfx))>0$ thus $M_1$ is invertible.
Hence using the determinant formula for block matrices and the definition of $w_{g,j}$ from \eqref{equ:def w, v} we obtain
\[
  \det(M) = \det(M_1) (M_3 - M_2^T M_1^{-1} M_2)
  = \det(M_1) \|w_{g,j}(\bfx)\|_2^2.
\]
Now, denoting $\lambda_1(M) \geq \cdots \geq \lambda_m(M) > 0$ and $\lambda_1(M_1) \geq \cdots \geq \lambda_{m-1}(M_1) > 0$ the eigenvalues of respectively $M$ and $M_1$, which are both symmetric positive semi-definite, the Cauchy interlacing theorem states that for all $1\leq k\leq m-1$ it holds $\lambda_k(M) \geq \lambda_k(M_1) \geq \lambda_{k+1}(M)$. 
We then have that 
\[
  \det(M) = \lambda_m(M) \prod_{k=1}^{m-1} \lambda_k(M)
  \geq \lambda_m(M) \prod_{k=1}^{m-1} \lambda_k(M_1)
  = \lambda_m(M) \det(M_1).
\]
Using the fact that $\det(M_1)>0$ (since $M_1$ is symmetric positive definite) and that $\lambda_m(M) = \sigma_m(\nabla g(\bfx))^2$, we obtain the desired inequality $ \sigma_m(\nabla g(\bfx))^2 \leq \|w_{g,j}(\bfx)\|_2^2$.
Finally, for $1\leq j\leq m$, consider $\tilde g \in \calC^1(\calX, \bbR^m)$ defined by $\tilde g_m = g_j$, $\tilde g_j = g_m$ and $\tilde g_i = g_i$ for $i\in \{1, \cdots, m\}\setminus\{j, m\}$. 
By invariance of the determinant under permutation, the previous inequality also holds for $j \neq m$.

\subsection{Proof of \Cref{prop:suboptimality uniform m}}
\label{subsec:proof prop:suboptimality uniform m}

Let $g\in\calG_m$ and $1\leq j\leq m$.
Firstly, let us state the large deviation inequality satisfied by $\|w_{g,j}(\bfX)\|_2^2$.
Since $\|w_{g,j}(\bfX)\|_2^2 \leq \|\nabla g_j(\bfX)\|_2^2$ and since $x\mapsto \|\nabla g_j(\bfx)\|_2^2$ is a nonzero polynomial of total degree at most $2\ell$, using \Cref{prop:bounded uniform quantiles} we have that 
\[
  \|\nabla g_j(\bfX)\|_2^2 \leq 2 \overline{\eta}_{A, s}^{2\ell} \Expe{\|\nabla g_j(\bfX)\|_2^2} \leq  2 \overline{\eta}_{A, s}^{2\ell} \kappa_g,
  \quad 
  \kappa_g := \sup_{1\leq j\leq m} \Expe{\|\nabla g_j(\bfX)\|_2^2},
\]
almost surely with $\kappa_g \leq 1$ by assumption on $\calG_m$.
Hence \eqref{equ:norm projection m} yields
\[
  \calL_{m,j}(g) 
  \leq \tilde{\gamma}_2
  \calJ(g),
  \quad
  \tilde{\gamma}_2 :=
  2 \overline{\eta}_{4, s}^{2\ell},
\]
where $\tilde \gamma_2 \leq 2^{1+6\ell} s^{-2\ell}$ since $\overline{\eta}_{4,s}=4(1-2^{-s})^{-1} \leq 8s^{-1}$.
Secondly, let us state the small deviation inequality satisfied by $\|w_{g,j}(\bfX)\|_2^2$.
Since $g$ is a non-constant polynomial of total degree at most $\ell + 1$, we have that $M : \bfx \mapsto \nabla g(\bfx)^T \nabla g(\bfx)$ is a non-zero polynomial of total degree at most $2\ell m$. 
Thus applying \Cref{prop:small deviation sigma m} to $M$, then using \Cref{lem:bounds on norm w}, then using $\sigma_m(\nabla g(\bfX)^T \nabla g(\bfX)^T) = \sigma_m(\nabla g(\bfX))^2$, we obtain
\[
\begin{aligned}
  &\Proba{\|w_{g,j}(\bfX)\|_2^2 \leq q_{\det(M)} \varepsilon}
  \leq \underline{\eta}_{M, s} \varepsilon^{1/2\ell m},
  \\
  &\underline{\eta}_{M, s} = 
  \overline{\eta}_{4,s}
  (2m^{-1})^{\frac{1}{4\ell}}
  \sup_{\bfx\in\calX} \|M(\bfx)\|_F^{\frac{m-1}{2 \ell m}}
\end{aligned}
\]
with $\underline{\eta}_{4,s}$, $\overline{\eta}_{4,s}$ and $q$ respectively defined in \Cref{prop:small deviation remez uniform}, \Cref{prop:large deviation remez uniform} and \eqref{equ:def quantile function}.
Then using $\|M(\bfX)\|_F \leq \|\nabla g(\bfX)\|_F^2$ and $\|\nabla g_j(\bfX)\|_2^2 \leq 2\overline{\eta}_{4,s}^{2\ell} \kappa_g \leq 2\overline{\eta}_{4,s}^{2\ell}$ almost surely, we have that $\|M(\bfX)\|_F \leq 2 m \overline{\eta}_{4,s}^{2\ell}$ almost surely.
Using $(m-1)/(1+2\ell m) \leq 1/2\ell$ we then obtain
\[
  \underline{\eta}_{M, s} 
  \leq 2
  \underline{\eta}_{4,s} 
  \overline{\eta}_{4,s}^{1-\frac{1}{m}}
  m^{\frac{1}{4\ell}}.
\]
Now, using the same reasoning as in \Cref{subsec:proof prop:bound J by L}, with $p=1$ since $\bfX$ has compact support and using the upper bound on $\underline{\eta}_{M,s}$, we obtain
\[
  \calJ(g) 
  \leq 
  2\underline{\eta}_{M,s}^{\frac{2\ell m}{1+2\ell m}}
  q_{\det(M)}^{-\frac{1}{1+2\ell m}}
  \calL_{m,j}(g)^{\frac{1}{1+2\ell m}}.
\]
Then using \Cref{prop:bounded uniform quantiles} with $p_1\geq 1$ we have $q_{\det(M)}^{-1} \leq \min\{\overline{\eta}_{4,s},6Ap_1\ell m\}^{2\ell m} \underline{\nu}_{\calK_m,p_1}^{-1}$ with $\calK_m$ defined in \eqref{equ:def Km}.
Then, combining the upper bounds on $\underline{\eta}_{M, s}$ and $q_{\det(M)}^{-1}$, we obtain
\[
  \calJ(g) \leq
  \tilde \gamma_1
  \underline{\nu}_{\calK_m,p_1}^{-\frac{1}{1+2\ell m}}
  \calL_{m,j}(g)^{\frac{1}{1+2\ell m}},
  \quad 
  \tilde \gamma_1 
  := 2
  (
    2 \underline{\eta}_{4,s} 
    \overline{\eta}^{1-\frac{1}{m}}
    m^{\frac{1}{4\ell}}
    \min\{\overline{\eta}_{A,s}, 6Ap_1\ell m\}
  )^{\frac{2\ell m}{1+2\ell m}},
\]
where $\tilde \gamma_1 \leq 2^9 s^{-1} m^{1/4\ell} \min\{s^{-1}, 3p_1 \ell m\}$ since $\underline{\eta}_{4,s} \overline{\eta}_{4,s} = A^2 s^{-1}$ and $\overline{\eta}_{A,s} \leq 2As^{-1}$.
Finally, let $g\in\calG_m$ and let $g^*$ be a minimizer of $\calL_{m,j}$ over $\calG_{m,g}^{(j)}$.
Then for all $h\in\calG_{m,g}^{(j)}$ we have
\[
  \calJ(g^*)
  \leq 
  \tilde{\gamma}_1 
  \big(
  \underline{\nu}_{\calK_m,p_1}^{-1}
  \calL_{m,j}(h)
  \big)^{\frac{1}{1+2\ell m}}
  \leq \tilde{\gamma}_3 
  \underline{\nu}_{\calK_m,p_1}^{-\frac{1}{1+2\ell m}}  
  \calJ(h)^{\frac{1}{1+2\ell m}}
\]  
where, using in particular $\underline{\eta}_{4,s} \overline{\eta}_{4,s}=4^2 s^{-1}$ and $\overline{\eta}_{4,s} \leq 8s^{-1}$, 
\[
\begin{aligned}
  \tilde{\gamma}_3 
  := \tilde{\gamma}_1 \tilde{\gamma}_2^{\frac{1}{1+2\ell m}}
  & = 2\big(
    2^{1+\frac{1}{2\ell m}} \underline{\eta}_{A,s} \overline{\eta}_{A,s} m^{\frac{1}{4\ell}} \min\{\overline{\eta}_{A,s}, 6A p_1 \ell m\}
    \big)^{\frac{2\ell m}{1+2 \ell m}}, \\
  & \leq 2^{10} m^{\frac{1}{4\ell}} s^{-1} \min\{s^{-1}, 3 p_1 \ell m\}.
\end{aligned}
  \]

\subsection{Proof of \Cref{prop:Lm is quadratic form}}
\label{subsec:proof prop:Lm is quadratic form}

Let $g\in\calG_m$ anf $h\in\calG_{m,g}^{(j)}$.
Since $v_{g,j}(\bfX) \in \spanv{\nabla g_{-j}(\bfX)}^{\perp}$, the projectors $\Piperp_{v_{g,j}(\bfX)}$ and $\Piperp_{\nabla g_{-j}(\bfX)}$ commute.
Hence $\|\Piperp_{v_{g,j}(\bfX)} w_{h,j}(\bfX) \|_2^2 =  \nabla h_j(\bfX)^T  \Piperp_{v_{g,j}(\bfX)} \Piperp_{\nabla g_{-j}(\bfX)} \nabla h_j(\bfX)$.
Also, $\nabla h_j(\bfX) = \nabla \Phi(\bfX) G_j$, thus
\[
  \calL_{m,j}(h) 
  = G_j^T \Expe{
    \|v_{g,j}(\bfX)\|_2^2
    \nabla \Phi(\bfX)^T \Piperp_{v_{g,j}(\bfX)} \Piperp_{\nabla g_{-j}(\bfX)} \nabla \Phi(\bfX)
  } G_j.
\]
Then since $\|v_{g,j}(\bfX)\|_2^2 \Piperp_{v_{g,j}(\bfX)} = \|v_{g,j}(\bfX)\|_2^2 I_d - v_{g,j}(\bfX) v_{g,j}(\bfX)^T$ we have
\[
  \calL_{m,j}(h) 
  = G_j^T \Expe{
    \nabla \Phi(\bfX)^T 
    \big(
      \|v_{g,j}(\bfX)\|^2 \Piperp_{\nabla g_{-j}(\bfX)} - v_{g,j}(\bfX) v_{g,j}(\bfX)^T
    \big)  
    \nabla \Phi(\bfX)
  } G_j,
\]
in other words $ \calL_{m,j}(h)  = G_j^T H_{g,j} G_j$ where $H_{g,j}\in\bbR^{K\times K}$ defined in \eqref{equ:def Hm matrix for Lm} is symmetric, as the sum of two symmetric matrices, and since $\calL_{m,j}(h) \geq 0$ by definition, it is positive semi-definite.

\end{appendices}

\section*{Acknowledgments}
The authors thank the anonymous reviewers for their valuable suggestions. 
This work is supported by funds from the Agence Nationale de la Recherche (ANR-21-CE46-0015).

\section*{Conflict-of-interest}
The authors declare no competing interests.

\section*{Funding} 
This project is funded by the ANR-DFG project COFNET (ANR-21-CE46-0015). This work was partially conducted within the France 2030 framework programme, Centre Henri Lebesgue ANR-11-LABX-0020-01. 

\bibliographystyle{plain}  
\bibliography{main}

\end{document}